\documentclass{article}

\usepackage[a4paper,margin=2.5cm]{geometry}

\usepackage{amsfonts}
\usepackage{amsmath}
\usepackage{amssymb}
\usepackage{amsthm}
\usepackage{authblk}
\usepackage{bbold}
\usepackage{circuitikz}
\usepackage{diagbox}
\usepackage{dirtytalk}
\usepackage{float}
\usepackage{graphicx}
\usepackage{makecell}
\usepackage{mathtools}
\usepackage{orcidlink}
\usepackage[labelfont=bf,textfont=it]{caption}
\usepackage{tabularx}
\usepackage{thmtools}
\usepackage{subfig}
\usepackage{tikz-cd}
\usepackage{xcolor}
\usepackage{xurl}
\usepackage{datetime}
\usepackage{newtxmath}
\usepackage{placeins}

\usepackage[dvipsnames]{xcolor}

\usepackage{tikz}
\usetikzlibrary{arrows.meta}
\usetikzlibrary{calc}
\usetikzlibrary{patterns}

\usepackage{hyperref}
\usepackage{url}
\hypersetup{
  colorlinks=true,
  urlcolor=blue,
  citecolor=blue,
  linkcolor=red
}

\declaretheoremstyle[
  spaceabove=3pt,
  headfont=\bfseries, 
  notefont=\bfseries,
  notebraces={(}{)}, 
  bodyfont=\normalfont, 
]{boldstyle}

\declaretheorem[name=Theorem,style=boldstyle]{Theorem}
\declaretheorem[name=Definition,style=boldstyle]{Definition}
\declaretheorem[name=Proposition,style=boldstyle]{Proposition}
\declaretheorem[name=Lemma,style=boldstyle]{Lemma}

\declaretheorem[name=Example,style=boldstyle]{Example}
\declaretheorem[name=Note,style=boldstyle]{Note}
\declaretheorem[name=Construction,style=boldstyle]{Construction}

\begin{document}

\newcommand{\term}[1]{\textit{\textbf{#1}}}

\newcommand{\N}{\mathbb{N}}
\newcommand{\RP}{\mathbb{RP}^2}

\newcommand{\A}{\mathcal{A}}
\newcommand{\C}{\mathcal{C}}
\newcommand{\D}{\mathcal{D}}

\renewcommand{\c}{\textbf{c}}
\newcommand{\cycles}[1]{\textbf{Cycles}(#1)}
\newcommand{\cycle}[1]{\c}
\newcommand{\cycleset}[1]{\{\cycle{#1}\}}
\renewcommand{\i}{\textbf{i}}
\newcommand{\indiff}{e}

\newcommand{\Prefs}[1]{\textbf{Prefs}(#1)}
\newcommand{\WeakPrefs}[1]{\textbf{WeakPrefs}(#1)}

\newcommand{\bigbb}[1]{\mathord{\scalebox{1.2}{\ensuremath{\vmathbb{#1}}}}}
\newcommand{\3}{\bigbb{3}}
\newcommand{\threeast}{\3^\ast}

\renewcommand{\P}[1]{\textbf{Weak}(#1)}
\newcommand{\PN}[1]{\P{#1}^N}
\newcommand{\PP}[1]{\textbf{Strict}(#1)}
\newcommand{\PPN}[1]{\PP{#1}^N}

\newcommand{\PPValid}[1]{\textbf{Valid}(#1)}
\newcommand{\PPCont}[1]{\textbf{Cont}(#1)}

\newcommand{\PCont}[1]{\textbf{WeakCont}(#1)}
\newcommand{\Surface}[1]{\mathcal{S}[#1]}


\newcommand{\U}{\mathcal{U}}
\newcommand{\V}{\mathcal{V}}
\newcommand{\W}{\mathcal{W}}

\newcommand{\Nerve}[1]{\mathcal{N}(#1)}
\newcommand{\NPP}[1]{\Nerve[\PP{Blah}]}
\newcommand{\NPPN}[1]{\Nerve[\PPN{Blah}]}

\newcommand{\ModNerve}{\mathcal{M}}
\newcommand{\NPPC}[1]{\ModNerve[\PPC{#1}]}


\newcommand{\DTop}[0]{D_\text{top}}
\newcommand{\DBot}[0]{D_\text{bot}}
\newcommand{\XCover}[0]{\widetilde{X}}
\newcommand{\SCover}[0]{\widetilde{S}}

\newcommand{\parsplit}{\vspace{6pt}}


\setlength{\parskip}{6pt}

\title{The Non-Orientable Topology of Condorcet’s Paradox}

\renewcommand\Affilfont{\small}

\author[1,2]{Ori Livson \orcidlink{0009-0001-8425-589X}\thanks{Corresponding author: ori.livson@sydney.edu.au}}
\author[3]{Siddharth Pritam\orcidlink{0000-0001-5673-0406}}
\author[1,2]{Mikhail Prokopenko\orcidlink{0000-0002-4215-0344}}

\affil[1]{The Centre for Complex Systems,  University of Sydney, NSW 2006, Australia}
\affil[2]{School of Computer Science, Faculty of Engineering, University of Sydney, NSW 2006, Australia}
\affil[3]{Faculty of Computer Science, Chennai Mathematical Institute, Siruseri, Tamil Nadu, 603103, India}

\date{}
\maketitle
\vspace{-24pt}

\begin{abstract}
  Preference cycles are prevalent in problems of decision-making, and are contradictory when preferences are assumed to be transitive. This contradiction underlies Condorcet's Paradox, a pioneering result of social choice theory, wherein intuitive and seemingly desirable constraints on decision-making necessarily lead to contradictory preference cycles. Topological methods have since broadened social choice theory and elucidated existing results. However, characterisations of preference cycles in topological social choice theory are lacking. In this paper, we address this gap by introducing a framework for topologically modelling preference cycles that generalises Baryshnikov's existing topological model of strict, ordinal preferences on 3 alternatives. In our framework, the contradiction underlying Condorcet's Paradox topologically corresponds to the non-orientability of a surface homeomorphic to either the Klein bottle or real projective plane, depending on how preference cycles are represented. These findings allow us to reformulate Arrow's Impossibility Theorem in terms of the orientability of a surface as well.
\end{abstract}

\maketitle

\section{Introduction}\label{section:introduction}

\subsection{Arrow's Impossibility Theorem and Condorcet's Paradox}

Arrow’s Impossibility Theorem is a seminal result of social choice theory that demonstrates the impossibility of ranked-choice decision-making processes (e.g., voting methods) to jointly satisfy a number of intuitive and seemingly desirable constraints~\cite{sep-arrows-theorem}. The result is widely considered to have pioneered the field of social choice theory~\cite{sen-dictionary}, which has proven valuable to studies of how social decision-making can be done rather than how it is done~\cite{possibility-of-social-choice}.


A significant motivation for Arrow's Impossibility Theorem originates in Condorcet's observation that the requirement of majority rule can lead social decision-making processes to produce \term{preference cycles} \cite{arrow-history}. Preference cycles refer to situations such that for alternatives $X$, $Y$, and $Z$, $X$ is strictly preferred to $Y$, $Y$ to $Z$, and $Z$ to $X$, which we denote as $X \prec Y \prec Z \prec X$. In the context of ranked-choice voting, preferences must be transitive, which implies that preference cycles are \term{contradictory}. This is because any and all strict preferences (e.g., $X \prec Y$ and $Y \prec X$) simultaneously follow in a transitive preference cycle on three alternatives, which contradicts the strictness of the relation $\prec$.

In fact, Arrow's Impossibility Theorem has been explicitly shown to be a generalisation of Condorcet’s Paradox, specifically by D'Antoni~\cite{dantoni} for strict preferences and by Livson and Prokopenko for the general case, i.e., allowing indifference between alternatives~\cite{paper0-arxiv}. In other words, these works establish that ranked-choice decision-making processes fail to jointly satisfy the constraints of Arrow's Impossibility Theorem lest they aggregate certain individuals' preferences to contradictory preference cycles.

\subsection{{Valid vs. Contradictory Preference Cycles}}

The usefulness of developing tools to better understand preference cycles extends beyond social choice theory. This is due to the prevalence of preference cycles in other domains, e.g., in the study of money pumps, Dutch books and the rationality of intransitive preferences~\cite{anand-intransitivity,money-pump-gustafsson,dutch-book-hajek}. Indeed, in many cases, intransitive preference cycles are \term{valid}, i.e., not contradictory. For example, consider the set of all possible rules one could ascribe to the well-known game ``{rock-paper-scissors}''. If for $X, Y \in \{\text{Rock}, \text{Paper}, \text{Scissors}\}$, we write $X \prec Y$ to denote that ``{$X$ beats $Y$}'', the usual rules for the game correspond to the following preference cycle.

\begin{equation}
  \text{Rock $\prec$ Scissors $\prec$ Paper $\prec$ Rock}\label{equation:rps}
\end{equation}

These preferences are clearly intransitive and hence the preference cycle is valid. For instance, despite the fact that rock beats scissors and that scissors beats paper, it does not follow that rock beats paper; to the contrary, paper beats rock. Not only does the preference cycle of Equation (\ref{equation:rps}) constitute an entirely valid set of rules for the game, the alternative rules defined by reversing the preferences of Equation (\ref{equation:rps}) constitute another valid preference cycle.
\begin{equation}
  \text{Rock $\prec$ Paper $\prec$ Scissors $\prec$ Rock}\label{equation:psr}
\end{equation}

Hence, when modelling preference cycles as we do so in this paper using topological methods, it is imperative to distinguish between preference cycles that are valid (i.e., intransitive) and preference cycles that are contradictory (i.e., under the assumption of transitivity).

\subsection{{Topological Social Choice Theory}}

Methods from topology have been leveraged to study various problems of social choice theory, comprising a field known as topological social choice theory, pioneered by Chichilnisky and Heal~\cite{chichilnisky-heal}. On the one hand, this is unsurprising because social choice theory does not only study inputs and outputs that vary discretely (e.g., electoral ballots) but also those that may vary continuously (e.g., utility functions). On the other hand, classical results of social choice theory can be recovered via topological methods, as  demonstrated by Baryshnikov~\cite{baryshnikov,baryshnikov-2,baryshnikov-gibbard}.

Baryshnikov's key insight is that a set of ranked-choice preferences can be modelled by a simplicial complex, a fundamental object of study in algebraic topology: specifically, simplices in a particular nerve complex. Although preference cycles are not explicitly modelled in this scheme (i.e., as simplices), Chia has observed that the boundaries of these nerve complex models correspond to preference cycles~\cite{topological-social-choice-summary} (p. 5). However, a limitation of this representation of preference cycles is that one cannot distinguish between models of preference relations with vs. without preference cycles, i.e., one cannot consider the simplicial complex with vs. without its boundaries,  nor can one distinguish between valid and contradictory preference cycles in this representation.

\subsection{{Topologically Modelling Preference Cycles}}

To address the limitation of simplicial complex-based models of ranked-choice preferences being unable to distinguish between sets of preference relations with vs. without preference cycles, let alone valid vs. contradictory preference cycles, we generalise Baryshnikov's topological model of ranked-choice preferences as follows.

Firstly, to distinguish between topological models of sets of preference relations with vs. without preference cycles, we introduce nerve complexes that include additional simplices that correspond to preference cycles. In this scheme, by including or excluding those simplices, we are able to distinguish between topological models of sets of preference relations with vs. without preference cycles, respectively. When preference cycles are explicitly represented by simplices in a nerve complex model, we say that preference cycles are \term{realised} in that model, and they are \term{unrealised} otherwise.

Secondly, we topologically model contradictory preference cycles distinctly from valid preference cycles as follows. Recalling that any and all strict preferences hold in a contradictory preference cycle on three alternatives, we argue that all contradictory preference cycles are then equivalent and should therefore be identified with respect to a single object. In our topological models of contradictory preference cycles, this identification corresponds to a gluing operation, which requires an orientation-reversing twist to account for the opposite direction that preference cycles may reference the alternatives (see Equations (\ref{equation:rps}) and (\ref{equation:psr})).

Therefore, a key contribution of this paper is the development of four different topological models of strict ranked-choice preferences and preference cycles over three alternatives: specifically, one model for each of the four possible combinations of preference cycles being either valid or contradictory and preference cycles being either realised or unrealised in the respective model. In our framework, each combination is modelled by a surface, uniquely up to homeomorphism (see Table~\ref{table:summary}).

\begin{table}[h!]
  \centering
  \resizebox{0.9\linewidth}{!}{
    \setcellgapes{2pt}
    \makegapedcells
    \begin{tabular}{c|cc}
      \textbf{Preference Cycles Property} & \textbf{Unrealised} & \textbf{Realised} \\ \hline
      \makecell{\textbf{Valid}                                                      \\(Intransitive)}                      & \makecell{Annulus / Cylinder $(S^1 \times [0,1])$                     \\(Proposition~\ref{proposition:unrealised-valid})}    & \makecell{Sphere $(S^2)$                                    \\(Proposition~\ref{proposition:realised-valid})}                      \\\hline
      \makecell{\textbf{Contradictory}                                              \\(Transitive)}              & \makecell{Klein bottle $(K)$                     \\(Theorem~\ref{theorem:unrealised-contradictory})}      & \makecell{real projective plane $(\RP)$                     \\(Theorem~\ref{theorem:realised-contradictory}) }
    \end{tabular}
  }
  \caption{Models of strict ranked-choice preferences and preference cycles over three alternatives, classified by whether preference cycles are valid or contradictory, and realised or unrealised.}\label{table:summary}
\end{table}

Importantly, among these four topological models of preference cycles, a model is non-orientable precisely when it models contradictory preference cycles. Specifically, the model with unrealised, contradictory preference cycles is homeomorphic to the Klein bottle, and our model with realised, contradictory preference cycles is homeomorphic to the real projective plane.

As an application of this framework, by using our topological models where preference cycles are contradictory and realised, we show that Arrow's Impossibility Theorem is equivalent to a statement about the orientability of a surface derivable from a Social Welfare Function (Theorems~\ref{theorem:arrow-reduction-1} and~\ref{theorem:arrow-reduction-2}). In other words, we provide a topological reformulation rather than an independent topological proof of Arrow's Impossibility Theorem.

This result addresses a conjecture of Chichilnisky that \say{the cyclical element of Condorcet's Paradox is indicative of a topological problem}~\cite{chichilnisky-space-of-preferences} (p. 165): the problem being non-orientability in these models of ranked-choice preferences and preference cycles.

Moreover, our classification of topological models of preference cycles may be applicable to the study of preference cycles in other domains, e.g., \textit{money pumps} and \textit{Dutch books}. Furthermore, there is an existing wide-ranging interest in the relationship between non-orientability and circular statements found in many well-known logical paradoxes. For example, the non-orientability of the Möbius strip has been used to informally represent logical paradoxes such as the Liar paradox, e.g., by Kauffman~\cite{kauffman-virtual-logic}, and it appears as an example of Hofstadter's concept of strange loops~\cite{geb-review,geb}.

\section{Background}

In this section, we provide  technical prerequisites in social choice theory and topology necessary for the constructions underlying our results in Section~\ref{section:results}, although the majority of the social choice theory background is only required for the reformulation of Arrow's Theorem in Section~\ref{section:arrow-non-orientability}. The remaining results (i.e., Sections \ref{section:strict-orders}--\ref{section:contradictory}) only require the social choice theory background of Sections~\ref{subsection:weak-and-strict-orders} and~\ref{subsection:reference-orientation} and the topology background of Section~\ref{section:topology}.

\subsection{Social Choice Theory}\label{section:background-arrow-classic}

\subsubsection{Weak and Strict Orders}\label{subsection:weak-and-strict-orders}

Arrow's Impossibility Theorem concerns the aggregation of \term{weak orders}, which are transitive and complete relations. A canonical example of which is a preferential voting ballot, wherein an individual (vote) is a ranking of alternatives from most to least preferred. Weak orders permit tied rankings (i.e., \term{indifference}) between alternatives. We use the term \term{strict order} to refer to a weak order without indifference.

\pagebreak

Given a fixed, finite set of alternatives $\A$, a weak order $R \subseteq \A \times \A$ can be represented by relation symbols $\prec, \sim$ and $\preceq$, where for every pair $a, b \in \A$, we write the following:
\begin{itemize}
  \item $a \sim_R b$ for \term{indifference} between $a$ and $b$, i.e., when $(a,b) \in R$ and $(b,a) \in R$.
  \item $a \prec_R b$ for $a$ being \term{strictly preferred} to $b$, i.e., when $(a,b) \in R$ and $(b,a) \notin R$.
  \item $a \preceq_R b$ for $a$ being \term{weakly preferred} to $b$, i.e., only assuming that $(a,b) \in R$ holds.
\end{itemize}

When $R$ is clear from context, we drop the subscript (i.e., simply write $\prec, \sim$ and $\preceq$) or may also write \say{$a \prec b$ in $R$}. The defining properties of a weak order are correspondingly the following:
\begin{description}
  \item[\textbf{Transitivity}] $\forall a,b,c \in \A$: If $a \preceq b$ and $b \preceq c$, then $a \preceq c$.
  \item[\textbf{Completeness}] $\forall a,b \in \A$: One of $a \prec b$, $b \prec a$ or $a \sim b$ holds.
\end{description}

Moreover, weak orders may be written as a permutation of $\A$'s elements interpolated by the $\prec$ or $\sim$ symbols: for example, if $\A = \{a,b,c\}$, $a \prec b \sim c$ denotes the weak order consisting of $a \prec b$, $b \sim c$ and $a \prec c$ (by transitivity). Strict orders are chains written entirely using the $\prec$ symbol.

\subsubsection{Arrow's Impossibility Theorem}\label{subsection:arrows-impossibility-theorem}

Given a fixed number $N \in \N$ of individuals, a \term{profile} is an $N$-tuple of weak orders. An example of a profile is an election, i.e., a tuple containing a single ballot from each individual. Note, each individual corresponds to a fixed index in the tuple across profiles. A \term{Social Welfare Function} is a function from a set of \term{valid} profiles to a single aggregate weak order (e.g., an election outcome). Invalid profiles are those that would otherwise fail to aggregate to a weak order, e.g., because they aggregate to a contradictory preference cycle.

\begin{Definition}\label{definition:fairness-conditions}
  A Social Welfare Function satisfies the following:
  \begin{itemize}
    \item \textbf{Unrestricted Domain}: When all profiles are valid with respect to it (i.e., can be aggregated).

    \item \text{\textbf{Unanimity}: When all individuals strictly preferring $a$ over $b$ implies that the aggregate does too.}

    \item \textbf{Independence of Irrelevant Alternatives (IIA)}: When the outcome of aggregation with respect to alternatives $a$ and $b$ only depends on the individual preferences with respect to $a$ and $b$.

    \item \textbf{Non-Dictatorship}: When there is no individual such that, irrespective of the profile, their strict preferences are always present in the aggregate outcome. If this condition fails, we say that the Social Welfare Function has a \term{Dictator}.
  \end{itemize}
\end{Definition}

\begin{Theorem}[Arrow's Impossibility Theorem]
  If a Social Welfare Function on at least three alternatives and two individuals satisfies the Unrestricted Domain, Unanimity and IIA, then it must have a Dictator.
\end{Theorem}

\noindent For examples of proofs of Arrow's Impossibility Theorem, see~\cite{geanakoplos,arrow-one-shot}.

\subsection{Condorcet's Paradox}\label{section:condorcet-paradox}\label{section:background-arrow-classic-2}

Condorcet's paradox refers to the phenomenon where voting methods on three or more alternatives cannot guarantee that winners are always preferred by a majority of voters. A canonical example of this is the observation that the method of pairwise majority voting aggregates the profile specified by Table~\ref{table:condorcet-plain} to a contradictory preference cycle. Pairwise majority voting is defined by ranking alternatives $x \prec y$ if more voters strictly prefer $x$ to $y$ than $y$ to $x$, and it is $x \sim y$ if there is a tie. If we apply this rule to the profile defined by Table~\ref{table:condorcet-plain}, we find that the majority of individuals strictly prefer $a$ to $b$ (individuals 1 and 3),  $b$ to $c$ (individuals 1 and 2), and $c$ to $a$ (individuals 2 and 3). Thus, aggregation yields a preference cycle $a \prec b \prec c \prec a$, which is contradictory because aggregated preferences are required to be transitive in an election outcome.

\begin{table}[h]
  \centering
  \begin{tabularx}{0.9\textwidth}{|p{1.8in}|X|X|X|}
    \hline
    \diagbox{Ranking}{Individual} & \textbf{1} & \textbf{2} & \textbf{3} \\\hline
    \textbf{1}                    & $a$        & $b$        & $c$        \\\hline
    \textbf{2}                    & $b$        & $c$        & $a$        \\\hline
    \textbf{3}                    & $c$        & $a$        & $b$        \\\hline
  \end{tabularx}
  \caption{A profile on 3 voters and 3 alternatives $\{a,b,c\}$ that pairwise majority voting aggregates to a preference cycle.}
  \label{table:condorcet-plain}
\end{table}

It is a straightforward exercise to show that pairwise majority voting satisfies Unanimity, IIA and Non-Dictatorship, but as we have seen, it may violate the Unrestricted Domain. In fact, it can be shown that all Social Welfare Functions satisfying those first three constraints necessarily violate the Unrestricted Domain due to the existence of profiles that aggregate to preference cycles. In other words, Arrow's Impossibility Theorem is formally a generalisation of Condorcet's Paradox. This was proven by D'Antoni in the special case, where all preferences are strict~\cite{dantoni}; Livson and Prokopenko later extended this result to prove Arrow's Impossibility Theorem in full as a generalisation of Condorcet's Paradox~\cite{paper0-arxiv}. In both cases, this is done by generalising the definition of a Social Welfare Function to functions that not only aggregate profiles to strict or weak orders but also possibly to preference cycles as well.

\subsection{Arrow's Theorem as a Generalisation of Condorcet's Paradox}\label{section:arrow-as-generalisation-of-condorcet}

In this section, we review the aforementioned recent work that shows Social Welfare Functions satisfying Unanimity, IIA and Non-Dictatorship necessarily aggregate some profiles to relations with preference cycles. This constitutes a failure of the Unrestricted Domain given the other constraints and hence a proof of Arrow's Impossibility Theorem. We begin by reviewing D'Antoni's original approach to proving this fact for strict orders and preference cycles in the three-alternative case, and then, we explain how this approach can be generalised to include additional alternatives and indifference between alternatives.

D'Antoni's approach leverages the following system of using binary-valued tuples to represent both strict orders and preference cycles. Indeed, for any strict order $\prec$ over $\3 = \{a_1, a_2, a_3\}$, we can represent $\prec$ by a tuple $(t_1,t_2,t_3)$ where each $t_i$ ranges over $\{0,1\}$ as follows:
\begin{align*}
  t_1 = 0 \iff a_1 \prec a_2 &  & t_2 = 0 \iff a_2 \prec a_3 &  & t_3 = 0 \iff a_3 \prec a_1
\end{align*}
and $t_i = 1$ for the reverse, i.e.,
\begin{align*}
  t_1 = 1 \iff a_2 \prec a_1 &  & t_2 = 1 \iff a_3 \prec a_2 &  & t_3 = 1 \iff a_1 \prec a_3
\end{align*}
\noindent There are $2^3 = 8$ possible binary three-tuples on $\3$, which includes each of the six possible strict orders on $\3$ as follows.
\begin{equation}\label{equation:strict-orders}
  \begin{array}{lll}
    (0,0,1) & a_1 \prec a_2 \prec a_3 & \qquad\qquad\qquad (0,1,0) \quad a_3 \prec a_1 \prec a_2  \\
    (0,1,1) & a_1 \prec a_3 \prec a_2 & \qquad\qquad\qquad (1,0,0) \quad a_2 \prec a_3 \prec a_1  \\
    (1,0,1) & a_2 \prec a_1 \prec a_3 & \qquad\qquad\qquad (1,1,0) \quad a_3 \prec a_2 \prec a_1.
  \end{array}
\end{equation}

\noindent The remaining tuples $(0, 0, 0)$ and $(1, 1, 1)$ of Equation~\ref{equation:strict-orders} represent the distinct intransitive preference cycles $a_1 \prec a_2 \prec a_3 \prec a_1$ and $a_1 \prec a_3 \prec a_2 \prec a_1$. In the next section, we discuss why these are the only two such preference cycles on $\3$.

\subsubsection{The Reference Orientation of Preference Cycles}\label{subsection:reference-orientation}

The two preference cycles $a_1 \prec a_2 \prec a_3 \prec a_1$ and $a_1 \prec a_3 \prec a_2 \prec a_1$ are the \textit{only} two possible preference cycles on $\3$ in the sense that \textit{continuing} a preference cycle on either side can produce only one of $a_1 \prec a_2 \prec a_3 \prec a_1$ or $a_1 \prec a_3 \prec a_2 \prec a_1$ lexicographically. For example, the preference cycle $a_2 \prec a_3 \prec a_1 \prec a_2$ can be \textit{continued} on the right side (in blue) to form the following:

\begin{equation}\label{equation:continuation}
  a_2\ \prec\ a_3 \prec a_1 \prec\ a_2\ {\color{blue}\prec\ a_3 \prec\ a_1}
\end{equation}

We find that $a_1 \prec a_2 \prec a_3 \prec a_1$ is produced on the right-hand side of Equation (\ref{equation:continuation}). Continuing the relation leftwards also produces $a_1 \prec a_2 \prec a_3 \prec a_1$, but importantly, neither continuation produces the \textit{other} preference cycle $a_1 \prec a_3 \prec a_2 \prec a_1$ lexicographically. We call these two representative preference cycles $a_1 \prec a_2 \prec a_3 \prec a_1$ and $a_1 \prec a_3 \prec a_2 \prec a_1$ of $\3$ the two \term{reference orientations} any intransitive preference cycle can have on $\3$.

When preference cycles are contradictory, i.e., under the assumption of transitivity, we can repeatedly apply transitivity to produce either preference cycle out of the other; this renders preference cycles non-orientable in the above sense and, as we shall see, in a topological sense as well.

\begin{Example}
  We begin with the preference cycle $a_1 \prec a_2 \prec a_3 \prec a_1$. We first conclude by transitivity that $a_1 \prec a_2$ and $a_2 \prec a_3$ holding in the original cycle implies $a_1 \prec a_3$. Likewise, $a_3 \prec a_1$ and $a_1 \prec a_2$ holding in the original cycle implies $a_3 \prec a_2$. Finally, $a_2 \prec a_3$ and $a_3 \prec a_1$ implies $a_2 \prec a_1$. Combining these three conclusions, we have $a_1 \prec a_3 \prec a_2 \prec a_1$, i.e., the aforementioned other cycle.
\end{Example}

We conclude this section by outlining another argument that all \textit{contradictory} preference cycles on three alternatives are equivalent. Recall that preference cycles such as $a_1 \prec a_2 \prec a_3 \prec a_1$ are contradictory under the assumption of transitivity due to them causing all strict preferences $a_n \prec a_m$ and $a_m \prec a_n$ to simultaneously hold. Hence, all contradictory preference cycles can be considered \textit{logically equivalent} in the same sense that all contradictions in classical logic are (e.g., the proposition \textit{false} or $C \wedge \neg C$). Livson and Prokopenko formalise this argument by identifying all contradictory preference cycles as the bottom element in a lattice of weak orders, ordered by the strictness of each order's underlying preferences~\cite[Section 3.1]{paper1-arxiv} analogous to $\textit{false}$ being the bottom element of any Lindenbaum Algebra of classical propositional logic.

Thus, identifying the two preference cycles of $\3$ as an equivalent object
\textbf{c}, we denote the set of strict orders and the contradictory preference cycle on $\3$ by the disjoint
union $\PPCont{\3} \coloneqq \PP{\3} \sqcup \{{\textbf{c}}\}$.

\subsubsection{Social Welfare Functions with Contradictory Cycles}\label{subsection:generalised-social-welfare-functions}

Let $\PP{\3}$ be defined as the set of strict orders on $\3$, $\PPN{\3}$ as the corresponding sets of profiles, and $\Prefs{\3}$ as the set of all binary-valued tuples of length 3 (i.e., representing strict orders and preference cycles). D'Antoni defines a Social Welfare Function on $\3$ as any function of the form $\PPN{\3} \rightarrow \Prefs{\3}$~\cite{dantoni}. We note that this definition assumes a precursor to the Unrestricted Domain, since the only way for aggregation to fail is for it to produce preference cycles. However, this assumption is not an issue for proving Arrow's Impossibility Theorem because it suffices to show that an ordinary Social Welfare Function satisfying the Unrestricted Domain (which is certainly included in D'Antoni's definition) cannot also satisfy Unanimity, IIA and Non-Dictatorship.

In~\cite[Section 3]{dantoni}, D'Antoni defines analogues of Arrow's conditions (Definition~\ref{definition:fairness-conditions}) for functions of the form $w: \PPN{\3} \rightarrow \Prefs{\3}$. For example, the Unrestricted Domain is the property that its image $im(w) = \PP{\3}$ contains no preference cycles. Moreover, Arrow's conditions are recovered when the Unrestricted Domain holds. Hence, Arrow's Impossibility Theorem is a consequence of the following.

\begin{Theorem}\label{theorem:arrow-refined-base-original}
  If a Social Welfare Function $w: \PPN{\3} \rightarrow \Prefs{\3}$ satisfies Unanimity, IIA and Non-Dictatorship, then $\exists p \in \PPN{\3}$ such that $w(p) = (0,0,0)$ or $w(p) = (1,1,1)$, i.e., the Unrestricted Domain fails.
\end{Theorem}

D'Antoni first proves the above in~\cite[Section 4.2]{dantoni}. To generalise the above theorem to hold for three or more alternatives, one can use functions of the form $\PPN{\A} \rightarrow \Prefs{\A}$, where $\Prefs{\A}$ consists of binary tuples of length $\binom{|\A|}{2}$, i.e., the number of all unordered pairs of elements in $\A$ without repeats. This reflects that for each pair of distinct alternatives $a,b \in \A$, we must record whether $a \prec b$ or whether $b \prec a$. $\Prefs{\A}$ includes strict orders and relations with preference cycles (that may or may not span all elements of $\A$). Indeed, for any triple of alternatives, IIA implies that $w$ is determined on that triple by a three-alternative function, we can thus apply Theorem~\ref{theorem:arrow-refined-base-original} to produce a preference cycle on that triple that violates the Unrestricted Domain regardless of $\A$. To incorporate indifference, we simply use ternary valued tuples of length $\binom{|\A|}{2}$, which we denote by $\WeakPrefs{\A}$ (see~\cite{paper0-arxiv} for further details).

However, recall that preference cycles are transitive and hence contradictory in the context of Arrow's Theorem. Hence, we identify all contradictory preference cycles (i.e., the tuples $(0,0,0)$ and $(1,1,1)$) by a single object \textbf{c}, which admits quotienting on $\Prefs{\3}$ that is bijective to the disjoint union $\PPCont{\3} \coloneqq \PP{\3} \sqcup \{{\textbf{c}}\}$ (see Section~\ref{subsection:reference-orientation}).

Furthermore, applying this quotienting process to the image of D'Antoni's Social Welfare Functions allows us to define a \term{Social Welfare Function with Contradictory Cycles (SWFC)} in the strict case on $\3$ as any function of the form $\PPN{\3} \rightarrow \PPCont{\3}$. This allows us to rephrase Theorem~\ref{theorem:arrow-refined-base-original} as the following.

\begin{Theorem}\label{theorem:arrow-refined-base}
  If an SWFC $w: \PPN{\3} \rightarrow \PPCont{\3}$ satisfies Unanimity, IIA and Non-Dictatorship, then $\exists p \in \PPN{\3}$ such that $w(p) = \textbf{c}$, i.e., the Unrestricted Domain fails.
\end{Theorem}

In Appendix~\ref{appendix:more-alternatives}, we handle indifference by instead quotienting $\WeakPrefs{\3}$ to $\PCont{\3} \coloneqq \P{\3} \sqcup\{{\textbf{c}}\}$.

\subsection{Preliminaries in Topological Social Choice Theory}\label{section:topology}

\subsubsection{General Topological Preliminaries}\label{section:topology-preliminaries}

In this section, we briefly recall the basic notions from topology that are used in topological social choice theory; for further details, see~\cite{hatcher2002algebraic}.

The first topological construction we outline is that of a simplicial complex, which is effectively a geometric object made out of vertices, edges, triangles, tetrahedra, etc. Simplicial complexes generalise the concept of a graph, which only encodes pairwise relationships. As such, simplicial complexes can be used to \textit{geometrically realise} combinatorial objects. This section culminates in Baryshnikov's geometric realisation of the set of strict orders on a fixed set of alternatives by a simplicial complex.

A \term{geometric $k$-simplex} is the convex hull of $k+1$ affinely independent points in $\mathbb{R}^d$, where $k$ is the \term{dimension} of the simplex and $d \geq k$. For example, a point is a $0$-simplex, an edge is a $1$-simplex, a triangle is a $2$-simplex, and a tetrahedron is a $3$-simplex. A \term{face} of a simplex is a subset of its points that is a simplex itself, e.g., each vertex and edge of a triangle is a face of that triangle. A \term{geometric simplicial complex} $K$ is a collection of geometric simplices that intersect only along their faces and is closed under the face relation, i.e., all faces of simplices in the collection are also part of that collection (see Figure~\ref{fig:simpcompexp} for an example). The usual topology of a geometric simplicial complex is the subset topology induced by the standard topology of $\mathbb{R}^d$.

A vertex-to-vertex map $\phi : K \rightarrow L$ between simplicial complexes $K$ and $L$ is called a \term{simplicial map} when it maps every simplex in $K$ to a simplex in $L$. For example, a simplicial map may map a triangle (2-simplex) in $K$ to a triangle, an edge (1-simplex), or a vertex (0-simplex) in $L$.

Moreover, the above definitions can be generalised for finite sets rather than points in $\mathbb{R}^d$. A family of finite sets $\Delta$ is called an \term{abstract simplicial complex} if for every set $X \in \Delta$ (called a \term{simplex}) and subset $Y \subseteq X$ (called a \term{face}), we have $Y \in \Delta$. In other words, an abstract simplicial complex is a family of sets closed under the subset relation. Abstract simplicial complexes are automatically closed under intersection as well, just as geometric simplicial complexes are. Moreover, through canonical simplicial maps, every abstract simplicial complex determines a geometric realisation that is unique up to homeomorphism. Thus, for the purposes of this paper, given an abstract simplicial complex $K$, we may simply refer to $K$ as a simplicial complex, i.e., referring to this geometric realisation.

\begin{figure}[H]
  \centering
  \begin{tikzpicture}[scale=0.85, transform shape]

    \filldraw[black] (0, 0) circle (2pt) node[anchor=east] {$v_1$};
    \filldraw[black] (2, 0) circle (2pt) node[anchor=north] {$v_2$};
    \filldraw[black] (1, 2) circle (2pt) node[anchor=south] {$v_3$};
    \filldraw[black] (3, 2) circle (2pt) node[anchor=west] {$v_3$};
    \filldraw[black] (4, 0) circle (2pt) node[anchor=west] {$v_5$};

    \draw[thick] (0,0) -- (2,0) node[midway, below] {};
    \draw[thick] (0,0) -- (1,2) node[midway, left] {};
    \draw[thick] (2,0) -- (1,2) node[midway, right] {};
    \draw[thick] (2,0) -- (4,0) node[midway, below] {};
    \draw[thick] (2,0) -- (3,2) node[midway, right] {};
    \draw[thick] (3,2) -- (4,0) node[midway, right] {};
    \draw[thick] (1,2) -- (3,2) node[midway, above] {};

    \filldraw[gray, opacity=0.3] (0,0) -- (2,0) -- (1,2) -- cycle;
    \filldraw[gray, opacity=0.3] (2,0) -- (4,0) -- (3,2) -- cycle;

  \end{tikzpicture}

  \caption{Example of a simplicial complex.}
  \label{fig:simpcompexp}
\end{figure}

In this paper, we are concerned with a particular simplicial complex known as a \term{nerve complex} of a family of sets. In essence, the nerve is a summary of the intersection pattern of a family of sets. When those sets are suitably chosen to cover a surface, certain properties of the nerve complex reflect properties of the original surface.

\pagebreak

Nerve complexes are defined as follows. Given a family of sets $\mathcal{C} = \{ C_i \}_{i \in I}$, the nerve complex of $\mathcal{C}$ is a simplicial complex $\Nerve{\C}$ with vertex set $I$, and for any finite subset $\{ i_0, \dots, i_k \} \subseteq I$ of vertices, we have a $k$-simplex if and only if $C_{i_0} \cap \dots \cap C_{i_k} \neq \emptyset$. Moreover, a family $\mathcal{C}$ \term{covers} a set $X$ if $X = \bigcup_{i \in I} C_i$. See Figure~\ref{fig:nerve-example} for examples of nerve complexes.

\begin{figure}[H]
  \centering
  \begin{tikzpicture}[scale=0.85, transform shape]
    \draw[blue,opacity=0.3,fill=blue!20] (0.5,0.5) circle (1);
    \node[blue] at (-0.5,1.5) {$U_1$};

    \draw[red,opacity=0.3,fill=red!20] (1.5,0.5) circle (1);
    \node[red] at (2.7,1.5) {$U_2$};

    \draw[green,opacity=0.3,fill=green!20] (1,1) circle (1);
    \node[green!50!black] at (1,2.5) {$U_3$};

    \filldraw[gray, opacity=0.3] (5,0) -- (7,0) -- (6,2) -- cycle;
    \filldraw[black] (5,0) circle (2pt) node[anchor=east] {$1$};
    \filldraw[black] (7,0) circle (2pt) node[anchor=west] {$2$};
    \filldraw[black] (6,2) circle (2pt) node[anchor=south] {$3$};
    \draw[thick] (5,0) -- (7,0);
    \draw[thick] (5,0) -- (6,2);
    \draw[thick] (7,0) -- (6,2);

    \def\yrow{-3}

    \draw[blue,fill=blue!20,opacity=0.3]
    plot[smooth cycle] coordinates {
        (0.0,\yrow) (0.8,\yrow+0.3) (0.6,\yrow+0.8) (0.0,\yrow+0.6) (-0.2,\yrow+0.3)
      };
    \node[blue] at (-0.4,\yrow+0.9) {$V_1$};

    \draw[red,fill=red!20,opacity=0.3]
    plot[smooth cycle] coordinates {
        (0.8,\yrow+0.2) (1.8,\yrow+0.4) (1.6,\yrow+0.9) (0.8,\yrow+0.7) (0.6,\yrow+0.4)
      };
    \node[red] at (2.0,\yrow+1.0) {$V_2$};

    \draw[green!50!black,fill=green!20,opacity=0.3]
    plot[smooth cycle] coordinates {
        (0.4,\yrow+0.7) (1.2,\yrow+1.0) (1.1,\yrow+1.5) (0.4,\yrow+1.3) (0.1,\yrow+1.0)
      };
    \node[green!50!black] at (0.9,\yrow+1.8) {$V_3$};

    \filldraw[black] (5,\yrow) circle (2pt) node[anchor=east] {$1$};
    \filldraw[black] (7,\yrow) circle (2pt) node[anchor=west] {$2$};
    \filldraw[black] (6,2+\yrow) circle (2pt) node[anchor=south] {$3$};
    \draw[thick] (5,\yrow) -- (7,\yrow);
    \draw[thick] (5,\yrow) -- (6,2+\yrow);

  \end{tikzpicture}

  \caption{Nerve complexes for families of sets $\{U_i\}_{i \in \{1,2,3\}}$ and $\{V_i\}_{i \in \{1,2,3\}}$ with different intersection patterns.}
  \label{fig:nerve-example}
\end{figure}

\subsubsection{Preferences as Simplicial Complexes}\label{section:preference-complexes}

The nerve construction relevant to social choice theory is made as follows. Let $\PP{\A}$ be the set of strict orders on alternatives $\A$. For any two alternatives $i, j \in \A$, we can define $U_{ij}$ as the set of strict orders where $i \prec j$ holds, i.e.,  $U_{ij} = \{p \in \PP{\A}\mid i \prec j \text{ in $p$ } \}$. Clearly, $\U = \{U_{ij}\}_{i,j \in \A}$ covers $\PP{\A}$. By taking intersections of the covering sets $U_{ij}$, we can increasingly specify orders in $\PP{\A}$. E.g., $U_{ij} \cap U_{jk}$ consists of all orders $p$ such that $i \prec j \prec k$ holds in $p$. Moreover, $U_{ij} \cap U_{jk} \cap U_{ki} = \emptyset$ because no preference cycle $i \prec j \prec k \prec i$ can hold in any $p \in \PP{\A}$.

The nerve complex $\Nerve{\U}$ for $\A = \{1,2,3\}$ is shown in Figure~\ref{fig:nerve-u}. We note that there is a 2-simplex precisely for every strict order over $\{1,2,3\}$, e.g., the 2-simplex with vertices $\{23,31,21\}$ is part of the complex because $U_{23} \cap U_{31} \cap U_{21} = \{2 \prec 3 \prec 1\}$, i.e., the intersection is non-empty. Moreover, $\Nerve{\U}$ has two boundaries given by vertices $\{13,32,21\}$ and $\{12,23,31\}$. They do not comprise 2-simplices because were the corresponding intersections of the covering sets non-empty, those 2-simplices would correspond to the preference cycle $1 \prec 3 \prec 2 \prec 1$ and $1 \prec 2 \prec 3 \prec 1$, respectively.

\begin{figure}[H]
  \centering
  \resizebox{0.3\textwidth}{!}{%
    \begin{circuitikz}
      \tikzstyle{every node}=[font=\normalsize]

      \coordinate (12) at (6,12.25);
      \coordinate (23) at (3,8);
      \coordinate (31) at (9,8);
      \coordinate (13) at (5.25,10);
      \coordinate (32) at (6.75,10);
      \coordinate (21) at (6,9);

      \node[above] at (12) {12};
      \node[below left] at (23) {23};
      \node[below right] at (31) {31};
      \node[above left] at (13) {13};
      \node[above right] at (32) {32};
      \node[below] at (21) {21};

      \draw[black, line width=1pt] (12) -- (23);
      \draw[black, line width=1pt] (23) -- (31);
      \draw[black, line width=1pt] (31) -- (12);
      \draw[black, line width=1pt] (13) -- (32);
      \draw[black, line width=1pt] (21) -- (13);
      \draw[black, line width=1pt] (21) -- (32);
      \draw[black, line width=1pt] (21) -- (23);
      \draw[black, line width=1pt] (13) -- (23);
      \draw[black, line width=1pt] (32) -- (31);
      \draw[black, line width=1pt] (31) -- (21);
      \draw[black, line width=1pt] (32) -- (12);
      \draw[black, line width=1pt] (13) -- (12);

      \fill[gray!40,opacity=0.6] (23) -- (13) -- (21) -- cycle;
      \fill[gray!40,opacity=0.6] (23) -- (13) -- (12) -- cycle;
      \fill[gray!40,opacity=0.6] (23) -- (21) -- (31) -- cycle;
      \fill[gray!40,opacity=0.6] (31) -- (21) -- (32) -- cycle;
      \fill[gray!40,opacity=0.6] (31) -- (32) -- (12) -- cycle;
      \fill[gray!40,opacity=0.6] (12) -- (13) -- (32) -- cycle;
    \end{circuitikz}
  }%

  \caption{The nerve complex $\Nerve{\U}$ of $\U = \{U_{ij}\}_{i,j \in \A}$ for the 3-alternative case $\A = \3 = \{1,2,3\}$; each of the Nerve's vertex labels $ij$ refers to an index of a covering set of $\U$.}
  \label{fig:nerve-u}
\end{figure}

Baryshnikov's key insight is that a Social Welfare Function $w: \PPN{\A} \rightarrow \PP{\A}$ satisfies IIA if and only if $w$ induces a simplicial map $\Nerve{\U}^N \rightarrow \Nerve{\U}$ for covering $\U = \{U_{ij}\}_{i,j \in \A}$  defined above. See~\cite{baryshnikov,baryshnikov-2,topological-social-choice-summary} for topological proofs of Arrow's Impossibility Theorem using this construction.

\pagebreak

\subsection{Non-Orientability}\label{section:non-orientability}

\begin{figure}[H]
  \centering
  \scalebox{0.45}{
    \begin{minipage}{\textwidth}
      \captionsetup[subfloat]{font=large,labelfont=large}
      \centering
      \subfloat[\centering]{
        \includegraphics[height=0.25\textheight,keepaspectratio]{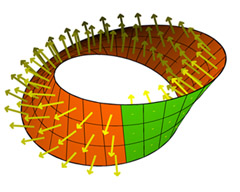}
      }
      \hspace{0.5em}
      \subfloat[\centering]{
        \includegraphics[height=0.25\textheight,keepaspectratio]{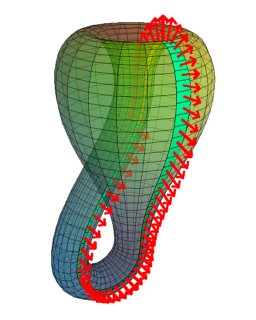}
      }

    \end{minipage}
  }

  \caption{The impossibility of orienting the M\"obius strip (\textbf{a}) and the Klein bottle (\textbf{b})~\cite{non-orientability-normals}.}
  \label{fig:non-orientability-normals}
\end{figure}

Orientation is defined by the assignment's normal vectors along a surface, and the opposing directions of these vectors correspond to opposing orientations, e.g., clockwise and anti-clockwise, top and bottom, and inside and outside. A surface is non-orientable when there is a path along a single orientation that leads to a loop that inverts the starting orientation (see Figure~\ref{fig:non-orientability-normals}). In this section, we outline how gluing operations (e.g., of nerve complexes) can be used to produce non-orientable surfaces.

\subsubsection*{Fundamental Polygons}\label{section:fundamental-polygons}

Examples of orientable surfaces include the annulus, sphere and the torus, while the M\"obius strip and the Klein bottle are well-known examples of non-orientable surfaces. These surfaces can be constructed from a square by stretching, twisting and gluing the sides of the square as illustrated in Figure~\ref{fig:fundamental-polygons} so that the labels and arrows.

\begin{figure}[H]
  \centering
  \resizebox{0.85\textwidth}{!}{%

    \begin{tikzpicture}[>=Stealth]


      \fill[gray!10] (0,0) rectangle (3,3);
      \draw[thick] (0,0) rectangle (3,3);

      \node at (-0.4,1.5) {$b$};
      \node at (3.5,1.5) {$b$};


      \draw[thick, ->] (-0.2,0) -- (-0.2,3);
      \draw[thick, ->] (3.2,0) -- (3.2,3);

      \node at (1.5,-1.1) {\textbf{Annulus/Cylinder}};

      \begin{scope}[shift={(5,0)}]

        \fill[gray!10] (0,0) rectangle (3,3);
        \draw[thick] (0,0) rectangle (3,3);

        \node at (-0.4,1.5) {$b$};
        \node at (3.5,1.5) {$b$};

        \draw[thick, ->] (-0.2,0) -- (-0.2,3);
        \draw[thick, <-] (3.2,0) -- (3.2,3);

        \node at (1.5,-1.1) {\textbf{M\"obius strip}};
      \end{scope}

      \begin{scope}[shift={(10,0)}]
        \fill[gray!10] (0,0) rectangle (3,3);
        \draw[thick] (0,0) rectangle (3,3);

        \node at (1.5,3.5) {$\alpha$};
        \node at (1.5,-0.4) {$\alpha$};
        \node at (-0.4,1.5) {$\beta$};
        \node at (3.5,1.5) {$\beta$};

        \draw[thick, ->] (0,3.2) -- (3,3.2);
        \draw[thick, ->] (0,-0.2) -- (3,-0.2);

        \draw[thick, <-] (-0.2,0) -- (-0.2,3);
        \draw[thick, ->] (3.2,0) -- (3.2,3);

        \node at (1.5,-1.1) {\textbf{Klein bottle}};
      \end{scope}

      \begin{scope}[shift={(15,0)}]
        \fill[gray!10] (0,0) rectangle (3,3);
        \draw[thick] (0,0) rectangle (3,3);

        \node at (1.5,3.5) {$\alpha$};
        \node at (1.5,-0.4) {$\alpha$};
        \node at (-0.4,1.5) {$\beta$};
        \node at (3.5,1.5) {$\beta$};

        \draw[thick, ->] (0,3.2) -- (3,3.2);
        \draw[thick, <-] (0,-0.2) -- (3,-0.2);

        \draw[thick, ->] (-0.2,0) -- (-0.2,3);
        \draw[thick, <-] (3.2,0) -- (3.2,3);

        \node at (1.5,-1.1) {\textbf{real projective plane}};
      \end{scope}

    \end{tikzpicture}
  }

  \caption{Fundamental polygons for various surfaces, wherein matching labels must be identified (i.e., glued) such that the arrow directions align.}

  \label{fig:fundamental-polygons}
\end{figure}

\begin{figure}[H]
  \centering
  \hspace{-3mm}\resizebox{0.85\textwidth}{!}{%
    \begin{tabular}{cccccc}
      \includegraphics[height=3cm]{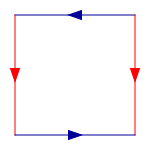} &
      \includegraphics[height=3cm]{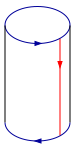} &
      \includegraphics[height=3cm]{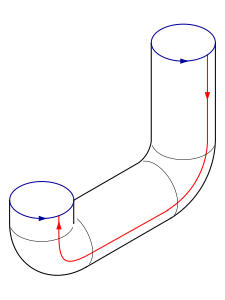} &
      \includegraphics[height=3cm]{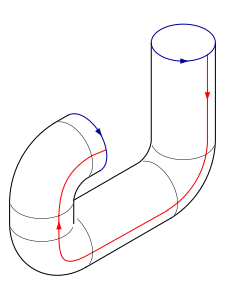} &
      \includegraphics[height=3cm]{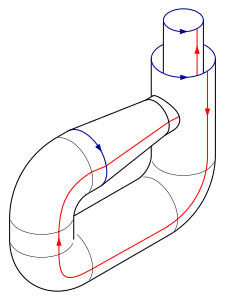} &
      \includegraphics[height=3cm]{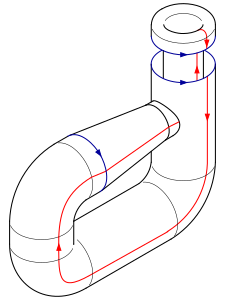}
    \end{tabular}%
  }

  \caption{Construction of a Klein bottle from its fundamental polygon~\cite{klein-construction-common}.}
  \label{fig:klein-construction}
\end{figure}

\begin{figure}[H]
  \centering
  \scalebox{0.45}{
    \begin{minipage}{\textwidth}
      \centering
      \captionsetup[subfloat]{font=large,labelfont=large}

      \subfloat[\centering]{
        \includegraphics[height=0.225\textheight,keepaspectratio]{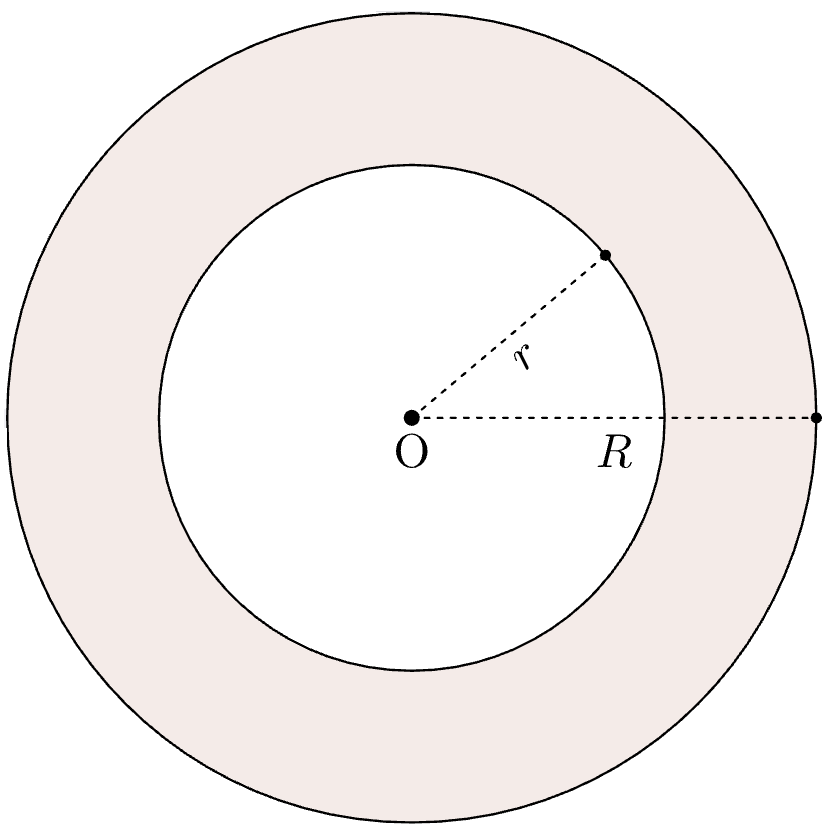}
      }
      \hspace{5em}
      \subfloat[\centering]{
        \includegraphics[height=0.225\textheight,keepaspectratio]{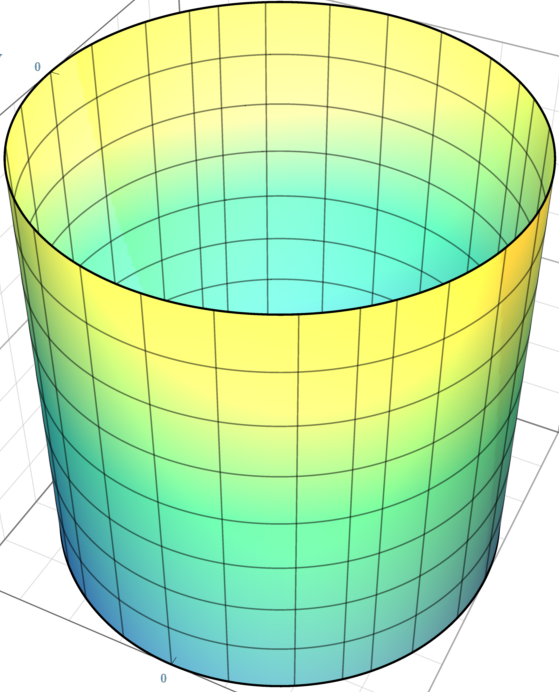}
      }

    \end{minipage}
  }
  \caption{Two surfaces homeomorphic to $S^1 \times [0,1]$, i.e., the product of a hollow circle and unit interval: the annulus (\textbf{a})~\cite{annulus} and hollow cylinder (\textbf{b})~\cite{cylinder}.}
  \label{fig:annulus}
\end{figure}

\begin{figure}[H]
  \centering
  \scalebox{0.45}{
    \begin{minipage}{\textwidth}
      \centering
      \captionsetup[subfloat]{font=large,labelfont=large}

      \subfloat[\centering]{
        \includegraphics[height=0.25\textheight,keepaspectratio]{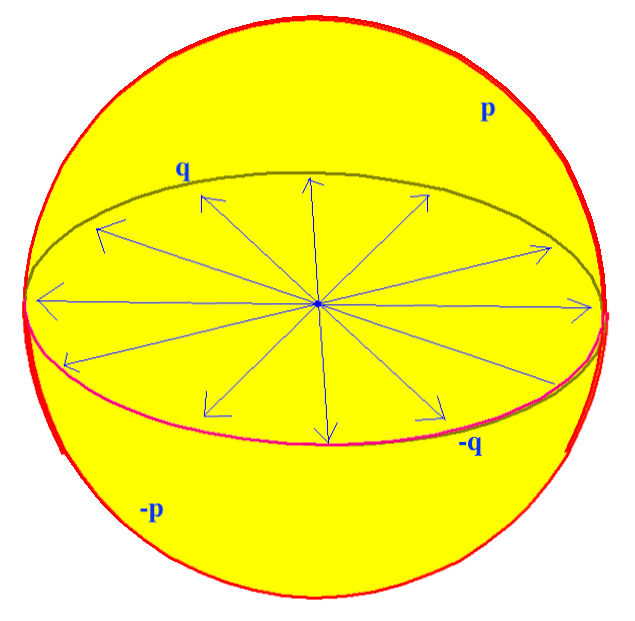}
      }
      \hspace{0.5em}
      \subfloat[\centering]{
        \includegraphics[height=0.25\textheight,keepaspectratio]{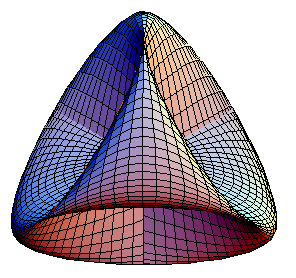}
      }

    \end{minipage}
  }
  \caption{(\textbf{a}) A construction of the real projective plane ($\RP$) by identifying antipodal points of a sphere~\cite{real-projective-plane-sphere}. (\textbf{b}) Visualisation of the Roman surface~\cite{roman-surface}, a surface homeomorphic to $\RP$.}
  \label{fig:real-projective-plane}
\end{figure}

The surfaces specified by the fundamental polygons of Figure~\ref{fig:fundamental-polygons}, can be constructed as follows. An (uncapped) cylinder can be constructed by rolling-up and gluing the two sides labelled $b$ together (the first step of Figure~\ref{fig:klein-construction} demonstrates this). Moreover, an annulus (see Figure~\ref{fig:annulus}) can be constructed by stretching a single side labelled $b$ out in an arc over the square until it meets the other side labelled $b$. A Möbius strip can be constructed in the same way, but one side needs to be twisted $180^{\circ}$ so that the arrows line up before gluing. The Klein bottle can be constructed per the steps of Figure~\ref{fig:klein-construction}, although those steps must be performed in 4D space to avoid \textit{self-intersection} of the surface. The construction of the real projective plane ($\RP$) can be done by identifying (i.e., gluing) antipodal points of a sphere (see Figure~\ref{fig:real-projective-plane}).

\begin{Note}\label{caveat-orientability}
  There is a separate notion of \term{orientable simplicial complexes}, which coincides with the usual notion of orientability with respect to triangulations of a surface. However, in this paper, we simply treat (geometric) simplicial complexes as surfaces and label their simplices to specify gluing operations as one would for a fundamental polygon.
\end{Note}

\section{Results}\label{section:results}

In Sections~\ref{section:strict-orders}--\ref{section:contradictory}, we define requirements and topological justifications for when a surface is a topological model of a set of strict orders and preference cycles on three alternatives. The requirements depend on whether the preference cycles being modelled are \term{valid} or \term{contradictory}; furthermore, preference cycles can either be \term{realised} or \term{unrealised} in each topological model (see Section~\ref{section:introduction}).

As we define these requirements, we identify solutions of topological models for each of the four combinations of preference cycles being valid or contradictory and realised or unrealised (see Table~\ref{table:summary}). Importantly, the solutions are non-orientable precisely for the topological models of contradictory preference cycles. In Section~\ref{section:arrow-non-orientability}, we apply this framework to reformulate Arrow's Impossibility Theorem in terms of the orientability of a surface derivable from a Social Welfare Function.

To begin, for a set $X$ containing strict orders and preference cycles on three alternatives, we will call a surface $\Surface{X}$ a \term{topological model of $X$} when it satisfies certain properties that depend on whether $X$ contains preference cycles and whether preference cycles are valid or contradictory, as well as realised or unrealised.

\subsection{Topological Models of Strict Orders}\label{section:strict-orders}

Let $\3 \coloneqq \{1,2,3\}$ and let $\PP{\3}$ be the set of all strict orders on $\3$. We begin by defining topological models of $\PP{\3}$. Recall Baryshnikov's definition of the following cover $\U = \{U_{ij}\}_{i,j\in\3}$ of $\PP{\3}$ with those covering sets defined as follows.
\begin{equation*}
  U_{ij} \coloneqq \{p \in \PP{\3} \mid i \prec j \text{ in $p$ } \}
\end{equation*}
Baryshnikov topologically represents $\PP{\3}$ by the nerve complex $\Nerve{\U}$, which has vertices denoted $\{ij\}_{i,j\in\3}$. Importantly, $\Nerve{\U}$ consists of six 2-simplices, one for each possible vertex set of the form $\{ab, bc, ac\}$, i.e., one 2-simplex per strict order $a \prec b \prec c$. Recall that the nerve complex $\Nerve{\U}$ has precisely those six 2-simplices because of the following identity (see Section~\ref{section:preference-complexes}).
\begin{equation*}
  U_{ab} \cap U_{bc} \cap U_{ac} = \{p \in \PP{\3} \mid a \prec b,\ b \prec c \text{ and } a \prec c \text{ in $p$ }\} = \{a \prec b \prec c\} \neq \emptyset
\end{equation*}

Since orientability is invariant under homeomorphism, we define topological models of $\PP{\3}$ up to homeomorphism as follows.

\begin{Definition}\label{definition:strict-order-top-mod}
  A surface $\Surface{\PP{\3}}$ is a topological model of $\PP{\3}$ when it is homeomorphic to $\Nerve{\U}$.
\end{Definition}

\parsplit

\begin{Proposition}\label{proposition:strict-order}
  If a surface $\Surface{\PP{\3}}$ is a topological model of $\PP{\3}$, then it is homeomorphic to the annulus/cylinder $S^1 \times [0,1]$, where $S^1$ denotes the circle (i.e., the 1D sphere).
\end{Proposition}
\begin{proof}
  By inspecting the visualisation of $\Nerve{\U}$ in Figure~\ref{fig:nerve-u}, we find that $\Nerve{\U} \cong S^1 \times [0,1]$. This is because one can stretch and round the edges of Figure~\ref{fig:nerve-u} to form an annulus or one can \textit{pull out} the inner boundary to form a cylinder.
  Specifically, $\Nerve{\U}$ can be given by $T \times [0,1]$ for an equilateral triangle $T$; the homeomorphism to the annulus can then be given by a transformation $(x_\theta, y_\theta) \mapsto (\frac{r}{\rho(\theta)}\cos\theta,\ \frac{r}{\rho(\theta)}\sin\theta)$, where $T$ is centred on $S^1$, $(x_\theta, y_\theta)$ is the point of the triangle $\theta$ degrees around $S^1$,  $\rho(\theta)$ is the distance from the centre of the triangle to $(x_\theta, y_\theta)$, and $r$ is the radius of $S^1$. Hence, for any topological model $\Surface{\PP{\3}}$ of $\PP{\3}$, we have $\Surface{\PP{\3}} \cong \Nerve{\U} \cong S^1 \times [0,1]$, as desired.\\
\end{proof}

In the next section, we generalise the above approach to topologically model a set containing strict orders on $\3$ and preference cycles on $\3$; firstly, under the assumption that preference cycles are \term{valid}. This means that the preference cycles $1 \prec 2 \prec 3 \prec 1$ and $1 \prec 3 \prec 2 \prec 1$ are distinct and permissible objects (see Section~\ref{section:introduction} for further examples).

\subsection{Topological Models of Strict Orders and Valid Preference Cycles}\label{section:valid}

Recall by Section \ref{subsection:reference-orientation} that there are only two distinct \term{valid} (i.e., intransitive) preference cycles on $\3$, which we denote as $\cycles{\3} \coloneqq \{1 \prec 2 \prec 3 \prec 1,\ 1 \prec 3 \prec 2 \prec 1\}$. Next, let the disjoint union $\PPValid{\3} \coloneqq \PP{\3} \sqcup \cycles{\3}$ denote the set of strict orders and valid preference cycles on $\3$.

Additionally, recall Chia's observation that the boundaries of $\Nerve{\U}$ given by vertex sets $\{12, 23, 31\}$ and $\{13, 32, 21\}$ correspond to the preference cycles $1 \prec 2 \prec 3 \prec 1$ and $1 \prec 3 \prec 2 \prec 1$, respectively~\cite{topological-social-choice-summary} (p. 5). Thus, surfaces homeomorphic to $\Nerve{\U}$ also model preference cycles in some manner. That is, $\Nerve{\U}$ consists of subspaces that model the elements of $\PP{\3}$ (i.e., its 2-simplices) and subspaces that model the elements of $\cycles{\3}$ (i.e., its two boundaries). However, $\Nerve{\U}$ is limited as a topological model of $\PPValid{\3}$ because it cannot be distinguished from a model of $\PP{\3}$, as its model of preference cycles (i.e., its boundaries) cannot be removed to produce a model of $\PP{\3}$ from $\PPValid{\3}$. We thus call preference cycles \term{unrealised} in the model $\Nerve{\U}$ of $\PPValid{\3}$.

Topological models of $\PPValid{\3}$ where preference cycles are \term{realised} will be defined as surfaces homeomorphic to a nerve complex $\Nerve{\V}$ for a covering $\V$ of $\PPValid{\3}$ rather than of $\PP{\3}$, which will consist of the simplices of $\Nerve{\U}$ and an additional 2-simplex per preference cycle in $\cycles{\3}$. The following covering $\V = \{V_{ij}\}_{i,j \in \3}$ of $\PPValid{\3}$ suffices.
\begin{equation}\label{equation:cover-2}
  \quad V_{ij} \coloneqq \{p \in \PPValid{\3} \mid i \prec j \text{ in $p$ and } j \nprec i \text{ in $p$}\}
\end{equation}

We proceed to verify that $\Nerve{\V}$ is equivalent to $\Nerve{\U}$ with the addition of the 2-simplices $\{12, 23, 31\}$ and $\{13, 32, 21\}$. Without loss of generality, for the strict order $1 \prec 2 \prec 3$, $1 \prec 2$ and $2 \prec 3$ hold, as well as $1 \prec 3$ by the transitivity property of strict orders. Moreover, we certainly have  $2 \nprec 1$, $3 \nprec 2$ and $3 \nprec 1$, which imply that $V_{12} \cap V_{23} \cap V_{13} = \{1 \prec 2 \prec 3\}$, as desired. Next, $\Nerve{\V}$ also has a 2-simplex for every valid preference cycle. Indeed, for a cycle $1 \prec 2 \prec 3 \prec 1$, we have that $1 \prec 2$, $2 \prec 3$ and $3 \prec 1$, but because in $\PPValid{\3}$, preference cycles are intransitive, we also have that $2 \nprec 1$ and $3 \nprec 2$ and also $3 \nprec 1$. Hence, $V_{12} \cap V_{23} \cap V_{31} = \{1 \prec 2 \prec 3 \prec 1\}$. This reasoning yields a 2-simplex for each of the eight elements of $\PPValid{\3}$, and furthermore, no other 2-simplices can exist in $\Nerve{\V}$ or else it would contain an edge of the form $\{ab, ba\}$ and clearly $V_{ab} \cap V_{ba} = \emptyset$.

\begin{Definition}\label{definition:valid-top-mod}
  A surface $\Surface{\PPValid{\3}}$ is a topological model of $\PPValid{\3}$ if and only if it is homeomorphic to either $\Nerve{\U}$ or $\Nerve{\V}$. When $\Surface{\PPValid{\3}}$ is homeomorphic to $\Nerve{\V}$, we say that the model has \term{realised} preference cycles and \term{unrealised} preference cycles otherwise.
\end{Definition}

\parsplit

\begin{Proposition}\label{proposition:unrealised-valid}
  If a surface $\Surface{\PPValid{\3}}$ is a topological model of $\PPValid{\3}$ with unrealised preference cycles, then it is homeomorphic to the annulus/cylinder $S^1 \times [0,1]$.
\end{Proposition}
\begin{proof}
  This follows identically to the proof of Proposition~\ref{proposition:strict-order}.\\
\end{proof}

\parsplit

\begin{figure}[H]
  \centering
  \includegraphics[width=0.45\textwidth]{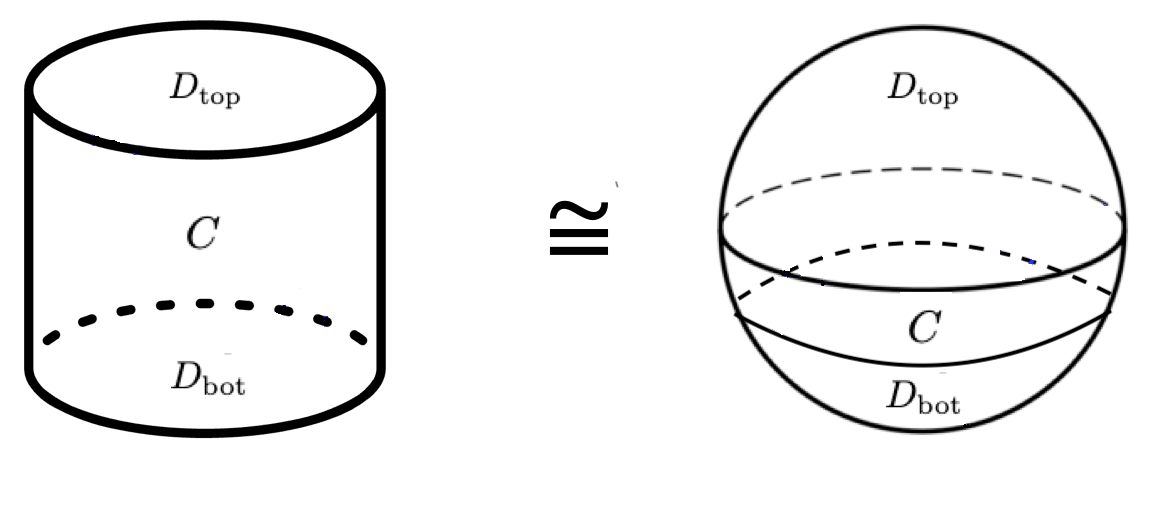}
  \caption{Homeomorphism of the capped cylinder $\overline{C} = \DTop \cup C \cup \DBot$ to the sphere $S^2$.}
  \label{fig:prop-sphere}
\end{figure}

\begin{Proposition}\label{proposition:realised-valid}
  If a surface $\Surface{\PPValid{\3}}$ is a topological model of $\PPValid{\3}$ with realised preference cycles, then it is homeomorphic to the 2D sphere, i.e., $S^2$.
\end{Proposition}
\begin{proof}
  Recall that, by definition, $\Surface{\PPValid{\3}} \cong \Nerve{\V}$ and that $\Nerve{\V}$ comprises $\Nerve{\U}$ with its two boundaries (i.e., vertex sets $\{12, 23, 31\}$ and $\{13, 32, 21\}$) filled in. Thus, by Proposition~\ref{proposition:unrealised-valid}, a topological model of $\Surface{\PPValid{\3}}$ is homeomorphic to the construction of attaching disks $\DTop$ and $\DBot$ to the two circular boundaries of the cylinder $C \coloneqq S^1 \times [0,1]$, {respectively,} resulting in the \textit{capped cylinder} $\overline{C} \coloneqq \DTop \cup C \cup \DBot$. Then, the homeomorphism from the capped cylinder to the sphere can be constructed by splitting the sphere into a top and bottom \say{hemisphere} with an \say{equatorial-band} between (see Figure~\ref{fig:prop-sphere}). The top hemisphere maps to $D_{top}$ by $(x,y,z) \mapsto (x,y,1)$, as does the bottom hemisphere to $D_{bot}$ with $(x,y,z) \mapsto (x,y,0)$. The remaining transformation of the equatorial band to the annulus $C$ can then be derived in a straightforward manner. Thus, we have that $\Surface{\PPValid{\3}} \cong \Nerve{\V} \cong S^2$.\\
\end{proof}

\subsection{Topological Models of Strict Orders and Contradictory Preference Cycles}\label{section:contradictory}

In this section, we begin by recalling the argument that, under the assumption of transitivity, the 
two distinct preference cycles of $\3$ are equivalent and should therefore be identified by a single object. Thus, to topologically model contradictory preference cycles, we begin with a model of valid preference cycles and identify (i.e., glue) the subspaces that correspond to those two preference cycles. However, we topologically demonstrate  that those subspaces must be assigned opposite orientations before identification, which in turn renders topological models of contradictory preference non-orientable.

To begin, recall that given transitivity, preference cycles on $\3$ are contradictory because it follows that any and all strict preferences (e.g., $1 \prec 3$ and $3 \prec 1$) must simultaneously hold. Thus, because the two preference cycles $\cycles{\3} = \{1 \prec 2 \prec 3 \prec 1, 1 \prec 3 \prec 2 \prec 1\}$ contain the same (i.e., all) strict preferences, they should be identified by a single object \textbf{c} (see Section~\ref{subsection:reference-orientation}). We denote the set of strict orders and contradictory preference cycles on $\3$ by the disjoint union $\PPCont{\3} \coloneqq \PP{\3} \sqcup \{{\textbf{c}}\}$, i.e., as a counterpart to the earlier defined $\PPValid{\3} \coloneqq \PP{\3} \sqcup \cycles{\3}$ for valid (i.e., intransitive) preference cycles. We then topologically model contradictory preference cycles by identifying the subspaces that represent each preference cycle in $\cycles{\3}$. However, we must first specify the orientations of those subspaces in order to identify them.

We first define criteria for when identifications of preference cycles in topological models of $\PPValid{\3}$ mirror the above logical identifications that produce $\PPCont{\3}$ as follows. We require that subspaces are only identified when they reflect an equivalence of preferences in $\3$, where preferences can be reversed by reversing a simplex's vertex order or reversing its assigned orientation. Specifically and without loss of generality, the 1-simplex $[12,23]$ with clockwise orientations denoting $12 \to 23$ can identify with either the counterclockwise-oriented $32 \to 21$ or equivalently with the clockwise-oriented $21 \to 32$. Moreover, identification of the preference cycle $12 \to 23 \to 31$ can be done in whole relative to any other preference cycle, e.g., $31 \to 12 \to 31$, with no restriction on orientation (since we seek to show that only the oppositely identified case is well-defined). However, this identification is indivisible because once $12 \to 23$ is identified, one can verify that the result can never connect to an identification of $23 \to 31$ regardless of the orientations used.

We next define criteria for how an identification of a pair of subspaces determines how its constituent (i.e., divided) parts are identified. It will then follow that the only consistent assignment of orientations to the subspaces representing the preference cycles of $\cycles{\3}$ is assigning them opposite orientations. The reason is that assigning them the same orientation will prevent the identification of parts that must be identified, while the same issue is not present when oppositely assigning orientations.

The setup for these criteria is as follows. First, recall that the preference cycles of $\cycles{\3}$ are topologically represented by the boundaries of $\Nerve{\U}$ or the 2-simplices with those boundaries in $\Nerve{\V}$. Thus, it suffices to demonstrate how the boundaries (i.e., preference cycles) of $\Nerve{\U}$ must identify because that determines how the interiors of the corresponding 2-simplices of $\Nerve{\V}$ would identify. Next, recall that $\Nerve{\U}$ is homeomorphic to the annulus $S^1 \times [0,1]$, and hence, $\Nerve{\U}$ can be represented by a fundamental polygon, as represented in Figure \ref{fig:nerve-fundamental-polygon}.

\begin{figure}[H]
  \centering
  \begin{tikzpicture}[scale=2.5]

    \coordinate (A) at (0,0);
    \coordinate (B) at (1,0);
    \coordinate (C) at (1,1);
    \coordinate (D) at (0,1);

    \draw (A) -- (B) -- (C) -- (D) -- cycle;

    \coordinate (P1) at (0,1);
    \coordinate (P2) at (0.333,1);
    \coordinate (P3) at (0.667,1);
    \coordinate (P4) at (1,1);

    \coordinate (Q1) at (0,0);
    \coordinate (Q2) at (0.333,0);
    \coordinate (Q3) at (0.667,0);
    \coordinate (Q4) at (1,0);

    \draw (Q1) -- (P1);
    \draw (Q4) -- (P4);

    \node[left]  at ($(Q1)!0.5!(P1)$) {$\alpha$};
    \node[right] at ($(Q4)!0.5!(P4)$) {$\alpha$};



    \foreach \pt in {P2,P3,Q1,Q2}{
        \fill (\pt) circle (0.01);
      }
    \foreach \pt in {P1,P4,Q1,Q4}{
        \fill (\pt) circle (0.02);
      }

    \node[above] at (P1) {13};
    \node[above] at (P2) {32};
    \node[above] at (P3) {21};
    \node[above] at (P4) {13};

    \node[below] at (Q1) {23};
    \node[below] at (Q2) {12};
    \node[below] at (Q3) {31};
    \node[below] at (Q4) {23};
  \end{tikzpicture}

  \caption{A fundamental polygon of $\Nerve{\U}$ with vertices at each corner. The edges labelled $\alpha$ are identified, and the top and bottom edges are currently the unidentified boundaries of $\Nerve{\U}$.}
  \label{fig:nerve-fundamental-polygon}
\end{figure}

The top and bottom sides of the polygon shown in Figure \ref{fig:nerve-fundamental-polygon} may be oriented and identified accordingly. For instance, we can assign the top side $[13,32,21,13]$ a clockwise orientation, and the counterclockwise orientation is assigned to the bottom side $[23,31,12,23]$. The orientation of any subspace modelling a preference cycle (i.e., a path with vertices $\{ab,bc,ca\}$) determines the orientation of the entire side to which it belongs. For instance, assigning an orientation to $[13,32,21]$ determines the full orientation of $[13,32,21,13]$. Furthermore, assigning an orientation to a cycle also determines the orientation of the remaining portion of the side it belongs to, e.g., the orientation $[13,32,21]$ determines the orientation $[21,13]$ of the top side of the polygon.


We can now prove that the subspaces representing preference cycles must be oppositely oriented before identification. Specifically, we first prove by contradiction that the top and bottom sides cannot be identified if they are assigned the same orientation; then, we show that assigning the sides opposite orientations is valid. Indeed, assume that the top and bottom sides have the same orientation: Without loss of generality, we say it is clockwise in Figure \ref{fig:nerve-fundamental-polygon}. This implies that we must identify $[13,32,21]$ and $[23,12,31]$, and in turn, we must identify $[21,13]$ and $[31,23]$ (see Figure \ref{fig:nerve-fundamental-polygon-2}). The pieces $[21,13]$ and $[31,23]$ encode the pairs of preferences $2 \prec 1 \prec 3$ and $2 \prec 3 \prec 1$, respectively. However, this identification violates our earlier criteria for topological identifications because $[21,13]$ and $[31,23]$ neither represent identical nor completely opposite preferences to one another.

We then prove that oppositely orienting the top and bottom sides of the polygon shown in Figure \ref{fig:nerve-fundamental-polygon} is valid by exhaustion. In other words, it suffices to verify that every orientation-reversing identification of the form $[ab,bc,ca]$ and $[ac,cb,ba]$ satisfies our earlier criteria for topological identifications. For instance, in Figure \ref{fig:nerve-fundamental-polygon-3}, the edges marked $\gamma$ both correspond to the two preference cycles with opposite reference orientations, and the edges marked $\delta$ correspond to $1 \prec 3 \prec 2$ on the top side and to the opposite $2 \prec 3 \prec 1$ on the bottom side, as desired. For a proof that all other cases are consistent, see Appendix~\ref{appendix:well-defined}.

\begin{figure}[H]
  \centering
  \begin{tikzpicture}[scale=2.5]
    \coordinate (P1) at (0,1);
    \coordinate (P2) at (0.333,1);
    \coordinate (P3) at (0.667,1);
    \coordinate (Q1) at (0,0);
    \coordinate (Q2) at (0.333,0);
    \coordinate (Q3) at (0.667,0);

    \coordinate (shift) at (0.4,0);
    \coordinate (P3r) at ($(P3)+(shift)$);
    \coordinate (P4) at ($(1,1)+(shift)$);
    \coordinate (Q3r) at ($(Q3)+(shift)$);
    \coordinate (Q4) at ($(1,0)+(shift)$);

    \draw (Q1) -- (P1) -- (P2) -- (P3) -- (Q3) -- (Q2) -- cycle;

    \draw (Q3r) -- (P3r) -- (P4) -- (Q4) -- cycle;

    \node[left] at ($(Q1)!0.5!(P1)$) {$\alpha$};
    \node[right] at ($(Q4)!0.5!(P4)$) {$\alpha$};

    \draw[thick] (P3) -- (Q3);
    \draw[thick] (P3r) -- (Q3r);
    \node[right] at ($(P3)!0.5!(Q3)$) {$\beta$};
    \node[left] at ($(P3r)!0.5!(Q3r)$) {$\beta$};

    \foreach
    \pt in {P2,Q2}{
        \fill (
        \pt) circle (0.01); }

    \foreach
    \pt in {P1,P3,P3r,P4,Q1,Q3,Q3r,Q4}{
        \fill (
        \pt) circle (0.02); }

    \node[above] at (P1) {13};
    \node[above] at (P2) {32};
    \node[above] at (P3) {21};
    \node[above] at (P3r) {21};
    \node[above] at (P4) {13};

    \node[below] at (Q1) {23};
    \node[below] at (Q2) {12};
    \node[below] at (Q3) {31};
    \node[below] at (Q3r) {31};
    \node[below] at (Q4) {23};

    \node[above=16pt] at ($(P1)!0.5!(P3)$) {$\gamma$};
    \draw[->] ($(P1)+(0,0.22)$) -- ($(P3)+(0,0.22)$);

    \node[below=16pt] at ($(Q1)!0.5!(Q3)$) {$\gamma$};
    \draw[->] ($(Q1)-(0,0.22)$) -- ($(Q3)-(0,0.22)$);

    \node[above=16pt] at ($(P3r)!0.5!(P4)$) {$\delta$};
    \draw[->] ($(P3r)+(0,0.22)$) -- ($(P4)+(0,0.22)$);

    \node[below=16pt] at ($(Q3r)!0.5!(Q4)$) {$\delta$};
    \draw[->] ($(Q3r)-(0,0.22)$) -- ($(Q4)-(0,0.22)$);

  \end{tikzpicture}

  \caption{A fundamental polygon equivalent to Figure \ref{fig:nerve-fundamental-polygon} with the top and bottom edges identified in an orientation-preserving manner. This is produced by subdividing \mbox{Figure \ref{fig:nerve-fundamental-polygon}} along $[21,31]$, which produces a new pair of identified edges $\beta$. The top and bottom sides of the polygon shown in \mbox{Figure \ref{fig:nerve-fundamental-polygon}} are identified with matching orientations if and only if the pairs $\gamma$ and $\delta$ are as well.}
  \label{fig:nerve-fundamental-polygon-2}
\end{figure}

\begin{figure}[H]
  \centering
  \begin{tikzpicture}[scale=2.5]
    \coordinate (P1) at (0,1);
    \coordinate (P2) at (0.333,1);
    \coordinate (P3) at (0.667,1);

    \coordinate (Q1) at (0,0);
    \coordinate (Q2) at (0.333,0);
    \coordinate (Q3) at (0.667,0);

    \coordinate (shift) at (0.6,0);

    \coordinate (P2r) at ($(P2)+(shift)$);
    \coordinate (P3r) at ($(P3)+(shift)$);
    \coordinate (P4)  at ($(1,1)+(shift)$);

    \coordinate (Q3r) at ($(Q3)+(shift)$);
    \coordinate (Q4)  at ($(1,0)+(shift)$);

    \draw (Q1) -- (P1) -- (P2) -- (Q3) -- (Q2) -- cycle;

    \draw (Q3r) -- (Q4) -- (P4) -- (P3r) -- (P2r) -- cycle;

    \node[left]  at ($(Q1)!0.5!(P1)$) {$\alpha$};
    \node[right] at ($(Q4)!0.5!(P4)$) {$\alpha$};

    \draw[thick] (P2) -- (Q3);
    \draw[thick] (P2r) -- (Q3r);

    \node[right]  at ($(P2)!0.5!(Q3)$) {$\beta$};
    \node[left] at ($(P2r)!0.5!(Q3r)$) {$\beta$};



    \foreach \pt in {P2,Q2}{
        \fill (\pt) circle (0.01);
      }
    \foreach \pt in {P1,P2,P2r,P4,Q1,Q3,Q3r,Q4}{
        \fill (\pt) circle (0.02);
      }

    \node[above] at (P1) {13};
    \node[above] at (P2) {32};
    \node[above] at (P2r) {32};
    \node[above] at (P3r) {21};
    \node[above] at (P4) {13};

    \node[below] at (Q1) {23};
    \node[below] at (Q2) {12};
    \node[below] at (Q3) {31};
    \node[below] at (Q3r) {31};
    \node[below] at (Q4) {23};

    \node[above=16pt] at ($(P1)!0.5!(P2)$) {$\gamma$};
    \draw[->] ($(P1)+(0,0.22)$) -- ($(P2)+(0,0.22)$);

    \node[below=16pt] at ($(Q1)!0.5!(Q3)$) {$\delta$};
    \draw[<-] ($(Q1)-(0,0.22)$) -- ($(Q3)-(0,0.22)$);

    \node[above=16pt] at ($(P2r)!0.5!(P4)$) {$\delta$};
    \draw[->] ($(P2r)+(0,0.22)$) -- ($(P4)+(0,0.22)$);

    \node[below=16pt] at ($(Q3r)!0.5!(Q4)$) {$\gamma$};
    \draw[<-] ($(Q3r)-(0,0.22)$) -- ($(Q4)-(0,0.22)$);
  \end{tikzpicture}

  \caption{A fundamental polygon equivalent to Figure \ref{fig:nerve-fundamental-polygon} with the top and bottom sides identified with opposite orientations. This is instead produced by subdividing the polygon along $[32,31]$. The top and bottom sides are identified if and only if the pairs $\gamma$ and $\delta$ are as well (aligning arrows).}
  \label{fig:nerve-fundamental-polygon-3}
\end{figure}
Hence, the following orientation of topological models $\Surface{\PPValid{\3}}$ is well-defined.

\parsplit

\begin{Definition}\label{definition:nerve-orientation}
  Let $\Surface{\PPValid{\3}}$ be a topological model of $\PPValid{\3}$ so that it is homeomorphic to a nerve complex with all 2-simplices and boundaries having vertex sets of the form $\{ab, bc, ac\}$ or $\{ab, bc, ca\}$, respectively. The \term{reference orientation} of any 2-simplex or boundary is the orientation that corresponds to the vertex ordering $[ab, bc, ac]$ or $[ab, bc, ca]$, respectively.
\end{Definition}

\begin{figure}[H]
  \centering

  \includegraphics[width=0.3\textwidth]{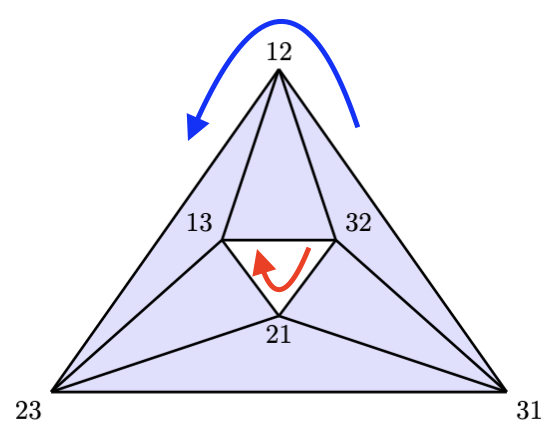}
  \caption{The reference orientation of $\Nerve{\U}$. 2-Simplex orientations are depicted by fill colour, and boundary orientations by arrow orientation and colour: in both cases, blue for counterclockwise and red for clockwise.}
  \label{fig:nerve-u-oriented}
\end{figure}

See Figure~\ref{fig:nerve-u-oriented} for an illustration of the reference orientation of $\Nerve{\U}$. To summarise, the reference orientation is an ordering of the vertices such that the alternatives appear lexicographically in the same order that they appear in the corresponding strict order or preference cycle. For instance, the strict order $a \prec b \prec c$ induces the lexicographic order $[ab, bc, ac]$ on the corresponding vertices. Moreover, a preference cycle $a \prec b \prec c \prec a$ having the reference orientation corresponding to $[ab, bc, ca]$ is well-defined because even though the cycle can be written as $b \prec c \prec a \prec b$ (with reference orientation corresponding to $[bc,ca,ab]$), the resulting orientation remains unchanged (i.e., both are clockwise or counterclockwise). See~Figure \ref{fig:nerve-arrow-realised} for an illustration of the reference orientation of $\Nerve{\V}$.

\noindent Importantly, in both $\Nerve{\U}$ and $\Nerve{\V}$, the two preference cycles are represented by boundaries or 2-simplices with opposite reference orientations. Moreover, the fact that all strict orders share an orientation is not problematic (with respect to our earlier criteria on identifications) because those 2-simplices are never identified.

\begin{figure}[H]
  \centering
  \resizebox{0.3\textwidth}{!}{%
    \begin{circuitikz}
      \tikzstyle{every node}=[font=\normalsize]

      \coordinate (12) at (6,12.25);
      \coordinate (23) at (3,8);
      \coordinate (31) at (9,8);
      \coordinate (13) at (5.25,10);
      \coordinate (32) at (6.75,10);
      \coordinate (21) at (6,9);

      \node[above] at (12) {12};
      \node[below left] at (23) {23};
      \node[below right] at (31) {31};
      \node[above left] at (13) {13};
      \node[above right] at (32) {32};
      \node[below] at (21) {21};

      \draw[black, line width=1pt] (12) -- (23);
      \draw[black, line width=1pt] (23) -- (31);
      \draw[black, line width=1pt] (31) -- (12);
      \draw[black, line width=1pt] (13) -- (32);
      \draw[black, line width=1pt] (21) -- (13);
      \draw[black, line width=1pt] (21) -- (32);
      \draw[black, line width=1pt] (21) -- (23);
      \draw[black, line width=1pt] (13) -- (23);
      \draw[black, line width=1pt] (32) -- (31);
      \draw[black, line width=1pt] (31) -- (21);
      \draw[black, line width=1pt] (32) -- (12);
      \draw[black, line width=1pt] (13) -- (12);

      \fill[blue!20,opacity=0.6] (23) -- (13) -- (21) -- cycle;
      \fill[blue!20,opacity=0.6] (23) -- (13) -- (12) -- cycle;
      \fill[blue!20,opacity=0.6] (23) -- (21) -- (31) -- cycle;
      \fill[blue!20,opacity=0.6] (31) -- (21) -- (32) -- cycle;
      \fill[blue!20,opacity=0.6] (31) -- (32) -- (12) -- cycle;
      \fill[blue!20,opacity=0.6] (12) -- (13) -- (32) -- cycle;

      \fill[pattern=north east lines, pattern color=red, line width=1pt] (13) -- (32) -- (21) -- cycle;
      \fill[pattern=north west lines, pattern color=blue, line width=1pt] (13) -- (32) -- (21) -- cycle;
    \end{circuitikz}
  }%

  \caption{The reference orientation of $\Nerve{\V}$ (see Equation~\ref{equation:cover-2}). 2-Simplex orientations are depicted by fill colour: blue for counterclockwise and red for clockwise, with blue-red cross-hatching to account for the overlap between the 2-simplices $[12,31,23]$ and $[13,32,21]$.}
  \label{fig:nerve-arrow-realised}
\end{figure}
Hence, we define topological models of $\PPCont{\3}$ as follows.

\parsplit

\begin{Definition}\label{definition:contradictory-top-mod}
  A surface $\Surface{\PPCont{\3}}$ is a topological model of $\PPCont{\3}$ when it is homeomorphic to the result of taking a topological model $\Surface{\PPValid{\3}}$ and identifying the subspaces that represent preference cycles according to their reference orientations. We say that preference cycles are \term{realised} in $\Surface{\PPCont{\3}}$ when they are realised in $\Surface{\PPValid{\3}}$ and \term{unrealised} otherwise.
\end{Definition}

We conclude this section by identifying solutions for topological models of $\PPCont{\3}$, i.e., depending on whether preference cycles are unrealised or realised in those models.

\parsplit

\begin{Theorem}\label{theorem:unrealised-contradictory}
  If a surface $\Surface{\PPCont{\3}}$ is a topological model of $\PPCont{\3}$ with unrealised preference cycles, then it is homeomorphic to the Klein bottle $K$.
\end{Theorem}
\begin{proof}
  By Definition~\ref{definition:contradictory-top-mod}, a topological model of $\PPCont{\3}$ with unrealised preference cycles is any surface that is homeomorphic to the result of beginning with $\Nerve{\U}$ and identifying its boundaries, i.e., the subspaces which correspond to the two preference cycles of $\cycles{\3}$. Moreover, this identification requires an orientation-reversing twist to align the opposite reference orientations of those two boundaries (see Figure~\ref{fig:nerve-u-oriented}).

  However, by Proposition~\ref{proposition:unrealised-valid}, this amounts to identifying the boundaries of a cylinder $S^1 \times [0,1]$ with an orientation-reversing twist, which comprises a well-known construction of the Klein bottle (see Figure~\ref{fig:klein-construction}).\\
\end{proof}

\begin{Theorem}\label{theorem:realised-contradictory}
  If a surface $\Surface{\PPCont{\3}}$ is a topological model of $\PPCont{\3}$ with realised preference cycles, then it is homeomorphic to the real projective plane $(\RP)$.
\end{Theorem}
\begin{proof}
  {By Definition~\ref{definition:contradictory-top-mod}, $\Surface{\PPCont{\3}}$ is any surface that is homeomorphic to the surface produced by taking $\Nerve{\V}$ and identifying its 2-simplices $\{12,23,31\}$ and $\{13,32,21\}$, i.e., the subspaces which correspond to the two preference cycles $1 \prec 2 \prec 3 \prec 1$ and $1 \prec 3 \prec 2 \prec 1$ of $\cycles{\3}$, respectively. Moreover, this identification requires an orientation-reversing twist to align the opposite reference orientations of those two 2-simplices (see Figure \ref{fig:nerve-arrow-realised}).

    Recall that in our proof of Proposition~\ref{proposition:realised-valid}, we showed that $\Nerve{\V}$ is homeomorphic to the capped cylinder $\overline{C} \coloneqq \DTop \cup C \cup \DBot$. Specifically, we derived a homeomorphism that associates the disks (i.e., the caps) $\DTop$ and $\DBot$ with the 2-simplices $\{12,23,31\}$ and $\{13,32,21\}$, respectively, and the remaining uncapped cylinder $C \subset \overline{C}$ with the remaining sub-simplicial complex $\Nerve{\U} \subset \Nerve{\V}$. Let $X$ denote the surface produced by taking $\overline{C}$ and identifying $\DTop$ and $\DBot$ with an orientation-reversing twist. Therefore, $\Surface{\PPCont{\3}} \cong X$.

    In Appendix~\ref{appendix:topology}, we prove that $X$ is homeomorphic to the real projective plane by explicitly parameterising $X$ in Euclidian space (see Construction \ref{construction}) and proving that $X \cong \RP$ (see Theorem \ref{theorem:main-construction}). This establishes that $\Surface{\PPCont{\3}} \cong \RP$.}\\
\end{proof}

We have thus identified all solutions for topological models for each of the four combinations of preference cycles being valid or contradictory, as well as realised or unrealised (see Table~\ref{table:summary}). We conclude Section~\ref{section:results} by applying these solutions to deriving a reformulation of Arrow's Impossibility Theorem in terms of the orientability of a surface derivable from a Social Welfare Function.

\subsection{Arrow's Impossibility Theorem and Non-Orientability}\label{section:arrow-non-orientability}

In this section, we reformulate Arrow's Impossibility Theorem in terms of the orientability of a surface derivable from a Social Welfare Function with Contradictory Cycles (SWFC) (see Section~\ref{subsection:generalised-social-welfare-functions}). We do so by beginning with the three-alternative case without indifference. The full results, i.e., with additional alternatives and indifference, are stated and proven in Appendix~\ref{appendix:more-alternatives}. We begin by summarising the result for the former case.

Recall that Arrow's Impossibility Theorem follows from the fact that an SWFC $w: \PPN{\3} \rightarrow \PPCont{\3}$ satisfying IIA, Unanimity and Non-Dictatorship necessarily aggregates some profile to a contradictory preference cycle (Theorem~\ref{theorem:arrow-refined-base}). In general, Unanimity implies that the image of an SWFC is either $\PP{\3}$ or $\PPCont{\3}$ (Lemma~\ref{lemma:two-options}). This allows us to define an \term{Arrovian Topological Model} of an SWFC $w$ as a topological model of $im(w)$ in the sense of Sections~\ref{section:strict-orders}--\ref{section:contradictory}: specifically, a topological model of $\PP{\3}$ (Definition~\ref{definition:strict-order-top-mod}) or a topological model of $\PPCont{\3}$ with realised preference cycles (\mbox{Definition~\ref{definition:contradictory-top-mod}}). Additionally, we say that $w$ has \term{(Non-)Orientable Social Choices} when its Arrovian Topological Models are (non-)orientable. Note, we define this for \say{models}, which is plural because the constituent topological models are surfaces defined only up to homeomorphism. This culminates in Theorem~\ref{theorem:arrow-reduction-1}, wherein Arrow's Impossibility Theorem is shown to be equivalent to the statement that no SWFC can simultaneously satisfy IIA, Unanimity, and Non-Dictatorship and have Orientable Social Choices.

We now proceed with the formalisation of the above argument. Fix a set $\3$ of three alternatives and $N \geq 2$. We begin with a lemma that allows us to define Arrovian Topological Models of an SWFC on $\3$ in a certain manner.

\parsplit

\begin{Lemma}\label{lemma:two-options}
  Let $w$ be an SWFC $w: \PPN{\3} \rightarrow \PPCont{\3}$ that satisfies Unanimity. Then, either $im(w) = \PP{\3}$ or $im(w) = \PPCont{\3}$.
\end{Lemma}
\begin{proof}
  By Unanimity, $\forall x \in \PP{\3}$: $w(x,x,\dots,x) = x$. This implies that $x \in im(w)$, which in turn implies that $\PP{\3} \subseteq im(w) \subseteq \PPCont{\3}$. Finally, because $\PPCont{\3} = \PP{\3} \sqcup\{{\textbf{c}}\}$ differs from $\PP{\3}$ by one element, there can be no proper subset between them; hence, $im(w) = \PP{\3}$ or $im(w) = \PPCont{\3}$.\\
\end{proof}

Note that the converse of Lemma \ref{lemma:two-options} does not hold knowing that the image of $w$ has no bearing on how $w$ behaves on profiles of the form $(x,x,\dots,x)$, as Unanimity requires. Using Lemma \ref{lemma:two-options}, we can define Arrovian Topological Models of $w$ as follows.

\begin{Definition}\label{definition:arrovian-topological-model}
  An \term{Arrovian Topological Model} of $w$ is any surface that is a topological model of $\PPCont{\3}$ with realised preference cycles when $im(w) = \PPCont{\3}$ and a topological model of $\PP{\3}$ otherwise, i.e., when $im(w) = \PP{\3}$.
\end{Definition}

Moreover, by the invariance of orientability under homeomorphism, we can additionally define the following property of SWFCs.

\begin{Definition}\label{definition:orientable-social-choices}
  Let $w: \PPN{\3} \rightarrow \PPCont{\3}$ be an SWFC on three alternatives. We say that $w$ has \term{Orientable Social Choices} when all of its Arrovian Topological Models are orientable and that $w$ has \term{Non-Orientable Social Choices} otherwise.
\end{Definition}

In fact, these two definitions are one and the same in the sense that all Arrovian Topological Models of $w$ are homeomorphic to the (orientable) annulus when $im(w) = \PP{\3}$ (Proposition~\ref{proposition:strict-order}) and homeomorphic to the (non-orientable) real projective plane ($\RP$) otherwise (Theorem~\ref{theorem:realised-contradictory}). We proceed to show that the two definitions are, in turn, one and the same as $w$ satisfying the Unrestricted Domain as well.

\begin{Proposition}\label{proposition:equivalent-conditions}
  An SWFC $w: \PPN{\3} \rightarrow \PPCont{\3}$ satisfies the Unrestricted Domain if and only if it has Orientable Social Choices.
\end{Proposition}
\begin{proof}
  ($\implies$) If $w$ satisfies the Unrestricted Domain, then $im(w) = \PP{\3}$, which implies that all Arrovian Topological Models of $w$ are topological models of $\PP{\3}$, which comprises any surface homeomorphic to the annulus (Proposition~\ref{proposition:strict-order}). Thus, because the annulus is orientable, we have that $w$ has Orientable Social Choices.

  ($\impliedby$) Assume by way of contradiction that $w$ has Orientable Social Choices but that $w$ does not satisfy the Unrestricted Domain. By Lemma~\ref{lemma:two-options}, if $w$ does not satisfy the Unrestricted Domain, then $im(w) = \PPCont{\3}$. This implies that all Arrovian Topological Models of $w$ must then be topological models of $\PPCont{\3}$ with realised preference cycles. However, this implies that they are homeomorphic to the (non-orientable) real projective plane (Theorem~\ref{theorem:realised-contradictory}), which contradicts our assumption that $w$ has Orientable Social Choices.\\
\end{proof}

\noindent Hence, we have the following equivalent reformulation of Arrow's Impossibility Theorem.

\parsplit

\begin{Theorem}[Reformulation of Arrow's Impossibility Theorem as Non-Orientability---Strict, 3-Alternative Case]\label{theorem:arrow-reduction-1}
  The following two statements are equivalent:
  \begin{enumerate}
    \item No SWFC $w: \PPN{\3} \rightarrow \PPCont{\3}$ satisfies all of Unanimity, IIA, Non-Dictatorship, and the Unrestricted Domain.
    \item No SWFC $w: \PPN{\3} \rightarrow \PPCont{\3}$ satisfies all of Unanimity, IIA, and Non-Dictatorship, and $w$ has Orientable Social Choices.
  \end{enumerate}
\end{Theorem}
\begin{proof}
  We simply apply Proposition~\ref{proposition:equivalent-conditions} to interchange the conditions of $w$ satisfying the Unrestricted Domain, and $w$ has Orientable Social Choices.\\
\end{proof}

This of course constitutes a reformulation of Arrow's Impossibility Theorem because statement 1 of Theorem~\ref{theorem:arrow-reduction-1} is exactly Arrow's Impossibility Theorem in the strict, three-alternative case (see Section~\ref{subsection:generalised-social-welfare-functions}). Arrow's Impossibility Theorem, i.e., for set $\A$ of three or more alternatives with indifference, is derived in Appendix~\ref{appendix:more-alternatives} in full; the approach can be summarised as follows.

Let $\3 \subseteq \A$ be any subset of three alternatives, and let $w$ be an SWFC on $\A$ for $N \geq 2$ individuals that satisfies IIA and Unanimity. IIA implies that the aggregate's preferences on $\3$ only depends on the individual's preferences on $\3$. Hence, $w$'s behaviour on $\3$ can be given by an SWFC $w_{\3}$ on just the alternative $\3$. We then define $\psi_{\3}$ by restricting the domain of $w_{\3}$ to the set of profiles $\D(\3)$ for which their aggregates are strict or have contradictory preference cycles, i.e., in $\PPCont{\3}$. This yields that $w$ satisfies the Unrestricted Domain if and only if $\psi_{\3}$ has Orientable Social Choices for every three-alternative subset $\3 \subseteq \A$ (Proposition~\ref{proposition:equivalent-conditions-2}). Then, we have a generalisation: Theorem~\ref{theorem:arrow-reduction-2} of Theorem~\ref{theorem:arrow-reduction-1} then follows in a straightforward manner.

To conclude, we prove the following stronger version of Arrow's Impossibility Theorem (Theorem \ref{theorem:stronger-equivalence-1}), which also has a corresponding generalisation in Appendix~\ref{appendix:more-alternatives}.

\parsplit

\begin{Lemma}\label{lemma:dictator-to-ud-base}
  If $w: \PPN{\3} \rightarrow \PPCont{\3}$ is an SWFC with a dictator, then it must satisfy the Unrestricted Domain.
\end{Lemma}
\begin{proof}
  Let $w$ be an SWFC with a dictator at individual $i$, let $p \in \PP{\3}$ be an arbitrary profile in the domain of $w$, and let $p_i \in \PP{\3}$ be the strict order comprising individual $i$'s preferences in $p$. By definition of $i$ being a dictator, $w(p)$ must contain every one of $i$'s strict preferences in $p$, which implies that $w(p) = p_i$. Hence, $im(w) = \PP{\3}$, which is the definition of $w$ satisfying the Unrestricted Domain.\\
\end{proof}

We note that the converse to Lemma \ref{lemma:dictator-to-ud-base} does not apply, as Non-Dictatorial rules may satisfy the Unrestricted Domain (e.g., instant run-off voting), but it may also fail the Unrestricted Domain (e.g., pairwise majority voting).

\parsplit

\begin{Theorem}[Non-Dictatorship as Non-Orientability (Strict, 3-Alternative Case)]\label{theorem:stronger-equivalence-1}~\\
  Let $w: \PPN{\3} \rightarrow \PPCont{\3}$ be an SWFC that satisfies IIA and Unanimity. The SWFC $w$ additionally satisfies Non-Dictatorship if and only if $w$ has Non-Orientable Social Choices.
\end{Theorem}
\begin{proof}
  ($\implies$) If $w$ satisfies Non-Dictatorship in addition to IIA and Unanimity, then, by Theorem~\ref{theorem:arrow-refined-base}, we have that $w$ does not satisfy the Unrestricted Domain, and hence, by Proposition~\ref{proposition:equivalent-conditions}, $w$ has Non-Orientable Social Choices.

  ($\impliedby$) If $w$ has Non-Orientable Social Choices, then by Proposition~\ref{proposition:equivalent-conditions} we have that $w$ does not satisfy the Unrestricted Domain, which in fact implies that $w$ cannot have a Dictator, i.e., by the contrapositive of Lemma~\ref{lemma:dictator-to-ud-base}.\\
\end{proof}

\noindent This result is generalised in full (for three or more alternatives with indifference) in Theorem~\ref{theorem:stronger-equivalence-2}.

\section{Discussion and Conclusions}

In this paper, we have established the topological problem of the cyclical element of Condorcet's Paradox---as Chichilnisky puts it~\cite{chichilnisky-space-of-preferences} (p. 165)---as the non-orientability of a topological model of ordinal preferences and preference cycles. To establish this result, we developed a framework for topological modelling of strict orders and preference cycles that generalises Baryshnikov's approach~\cite{baryshnikov,baryshnikov-2}. Within this framework, the topological properties of the corresponding surface depended on assumptions regarding preference cycles. In particular, when preference cycles were \term{contradictory} due to the assumption that preferences are transitive (as is the case in Condorcet’s Paradox), the resulting topological models were shown to be non-orientable. In contrast, when preference cycles were \term{valid} (as is the case in intransitive games such as rock–paper–scissors), the resulting topological models were shown to be orientable.

The fundamental distinction between contradictory and valid preference cycles arose out of the argument that all contradictory preference cycles on three alternatives are logically equivalent and should therefore be identified with a single object (see Section~\ref{subsection:reference-orientation}). Topologically, this identification of preference cycles corresponded to the identification of subspaces in a topological model, with an orientation-reversing twist. The orientation-reversing twist accounts for the opposite direction that alternatives are lexicographically referenced in each of the two possible preference cycles, i.e., $1 \prec 2 \prec 3 \prec 1$ vs. $3 \prec 2 \prec 1 \prec 3$. No such identification applies to valid preference cycles.

Another consequential property of these topological models was whether preference cycles corresponded to boundaries of a surface or to patches within the surface. In the former case, preference cycles were said to be \term{unrealised} in the model; in the latter, they were said to be \term{realised} in the model. When contradictory preference cycles were unrealised, the resulting non-orientable surface was found to be homeomorphic to the Klein bottle (Theorem~\ref{theorem:unrealised-contradictory}). Otherwise, when contradictory preference cycles were realised, the resulting non-orientable surface was shown to be homeomorphic to the real projective plane (Theorem~\ref{theorem:realised-contradictory}). The introduction of the concept of realised preference cycles was essential for distinguishing between topological models of sets of preference relations with vs. without preference cycles. This distinction was critical for reformulating Arrow's Impossibility Theorem (which generalises Condorcet's Paradox~\cite{dantoni,paper0-arxiv}) in terms of the orientability of a surface derivable from a Social Welfare Function (Theorem~\ref{theorem:arrow-reduction-1}).

\subsection{Related Work}
Several studies have developed proofs of Arrow's Impossibility Theorem that identify preference cycles as the boundaries of Baryshnikov's Nerve Complex representation of strict orders. This includes works by Chia~\cite{topological-social-choice-summary} and notably Raventós-Pujol and Rajsbaum~\cite{distributed-combinatorial-topology}, who use additional techniques from combinatorial topology and a related notion of orientable simplicial complexes (see Note~\ref{caveat-orientability}). However, these works do not directly relate these preference cycles to Condorcet's Paradox. In other words, they do not show that Social Welfare Functions satisfying Unanimity, IIA and Non-Dictatorship necessarily map some profile to these boundaries, nor do they make the distinction between valid and contradictory preference cycles that leads to modelling preference cycles in terms of non-orientability.

Ghrist has hypothesised that a different sort of topological impossibility underlies contradictory preference cycles, namely a cohomological obstruction~\cite{ghrist-textbook} (p. 39). Specifically, this is in reference to Penrose's cohomological explanation of the impossibility of constructing a Penrose triangle or the well-known impossible staircase~\cite{penrose, torsor}. The impossibility of the loop of ascending stairs can intuitively be likened to the impossibility of the loop of ascending preferences. However, this approach has not yet been applied to the concept of aggregation in the sense of Condorcet's Paradox.

Finally, Candeal and Indurian have previously reduced a special case of Chichilnisky's Impossibility Theorem to the non-retractability of the Möbius strip to its own boundary~\cite{candeal-indurian}. However, the result does not concern ranked-choice preferences and, as such, has not been related to preference cycles.

\subsection{Further Research Directions}
We conclude by discussing a number of promising topics for further research. Firstly, other results regarding preference cycles may be amenable to a topological characterisation as we have provided for Arrow's Impossibility Theorem. For example, preference cycles are prevalent in other domains~\cite{anand-intransitivity}, e.g., money pumps and Dutch books~\cite{may-intransitvity,money-pump-gustafsson,dutch-book-hajek}.

Secondly, visual paradoxes such as the Penrose triangle have been studied using topological methods~\cite{torsor,penrose}. As with Condorcet's Paradox, these paradoxes often arise out of issues regarding the assumption of transitivity. Hence, the topological methods used to study visual paradoxes may be applicable to the study of topologically modelling preference cycles and vice versa.

Lastly, the results of this paper suggest broader connections of non-orientability to circularity and self-reference, as found in the study of well-known logical paradoxes. For instance, recall that in this paper, non-orientability arose from the identification of all contradictory preference cycles as one and the same. In fact, this argument was previously leveraged by Livson and Prokopenko to establish formal relationships between Gödel's First Incompleteness Theorem and Arrow's Impossibility Theorem~\cite{paper1-arxiv}. In their analysis, Gödel sentences in logic play the same role as profiles that aggregate to contradictory preference cycles. The resulting triad of relationships between non-orientability, contradictory preference cycles and Gödel sentences (which are tantamount to the Liar paradox~\cite{prokopenko-liar-paradox}) suggests the existence of more fundamental relationships underlying circularity, self-reference and incomputability. Formal relations of which have been hypothesised by many other authors (e.g.,~\cite{tangled-hierarchies,kauffman-virtual-logic,geb}).

\FloatBarrier
\appendix
\setcounter{Theorem}{0}
\setcounter{Definition}{0}
\setcounter{Proposition}{0}
\setcounter{Lemma}{0}
\setcounter{Corollary}{0}
\setcounter{Example}{0}
\setcounter{Note}{0}
\setcounter{Construction}{0}

\renewcommand{\theTheorem}{\thesection.\arabic{Theorem}}
\renewcommand{\theDefinition}{\thesection.\arabic{Definition}}
\renewcommand{\theProposition}{\thesection.\arabic{Proposition}}
\renewcommand{\theLemma}{\thesection.\arabic{Lemma}}
\renewcommand{\theCorollary}{\thesection.\arabic{Corollary}}
\renewcommand{\theExample}{\thesection.\arabic{Example}}
\renewcommand{\theNote}{\thesection.\arabic{Note}}
\renewcommand{\theConstruction}{\thesection.\arabic{Construction}}

\renewcommand{\theHTheorem}{\thesection.\arabic{Theorem}}
\renewcommand{\theHDefinition}{\thesection.\arabic{Definition}}
\renewcommand{\theHProposition}{\thesection.\arabic{Proposition}}
\renewcommand{\theHLemma}{\thesection.\arabic{Lemma}}
\renewcommand{\theHCorollary}{\thesection.\arabic{Corollary}}
\renewcommand{\theHExample}{\thesection.\arabic{Example}}
\renewcommand{\theHNote}{\thesection.\arabic{Note}}
\renewcommand{\theHConstruction}{\thesection.\arabic{Construction}}

\section{Appendix}\label{appendix}

\subsection[\appendixname~\thesection]{Well-Definedness of Oppositely Oriented Preference Cycles}\phantomsection\label{appendix:well-defined}

In this section, we complete the argument of Section \ref{section:contradictory} that subspaces corresponding to preference cycles on $\3$ must be assigned opposite orientations before identification. So far, we have proven that those subspaces cannot be assigned the same orientation or else certain constituent identifications fail to be possible. It remains to show that oppositely orienting the subspaces is itself a valid possibility.

In order to prove this, we verify that every valid subdivision of the polygon shown in Figure \ref{fig:nerve-fundamental-polygon} leaving an orientation-reversing identification of cycles has all required identifications. In Section \ref{section:contradictory}, we proved one such case, i.e., if $[32,21,13]$ is identified with $[31,12,23]$, then $[13,32]$ is correctly identified with $[23,31]$. There are two other cases to check, visualised by Figures \ref{fig:nerve-fundamental-polygon-4} and \ref{fig:nerve-fundamental-polygon-5}. Indeed, in Figure \ref{fig:nerve-fundamental-polygon-4}, we have that the oppositely oriented edges marked $\gamma$ correspond to oppositely oriented cycles, and the edges marked $\delta$ correspond to $3 \prec 2 \prec 1$ on the top side and the reversed $1 \prec 2 \prec 3$ on the bottom side, as desired. Likewise, in Figure \ref{fig:nerve-fundamental-polygon-5}, we have that the oppositely oriented edges marked $\gamma$ correspond to oppositely oriented cycles, and the edges marked $\delta$ correspond to $2 \prec 1 \prec 3$ on the top side and to the reversed $3 \prec 1 \prec 2$ on the bottom side, as desired.

\begin{figure}[H]
  \centering
  \begin{tikzpicture}[scale=2.5]
    \coordinate (P1) at (0,1);
    \coordinate (P2) at (0.333,1);
    \coordinate (P3) at (0.667,1);

    \coordinate (Q1) at (0,0);
    \coordinate (Q2) at (0.333,0);
    \coordinate (Q3) at (0.667,0);

    \coordinate (shift) at (0.6,0);

    \coordinate (P2r) at ($(P2)+(shift)$);
    \coordinate (P3r) at ($(P3)+(shift)$);
    \coordinate (P4)  at ($(1,1)+(shift)$);

    \coordinate (Q3r) at ($(Q3)+(shift)$);
    \coordinate (Q4)  at ($(1,0)+(shift)$);

    \draw (Q1) -- (P1) -- (P2) -- (Q3) -- (Q2) -- cycle;

    \draw (Q3r) -- (Q4) -- (P4) -- (P3r) -- (P2r) -- cycle;

    \node[left]  at ($(Q1)!0.5!(P1)$) {$\alpha$};
    \node[right] at ($(Q4)!0.5!(P4)$) {$\alpha$};

    \draw[thick] (P2) -- (Q3);
    \draw[thick] (P2r) -- (Q3r);

    \node[right]  at ($(P2)!0.5!(Q3)$) {$\beta$};
    \node[left] at ($(P2r)!0.5!(Q3r)$) {$\beta$};



    \foreach \pt in {P3r,Q2}{
        \fill (\pt) circle (0.01);
      }

    \foreach \pt in {P1,P2,P2r,P4,Q1,Q3,Q3r,Q4}{
        \fill (\pt) circle (0.02);
      }

    \node[above] at (P1) {32};
    \node[above] at (P2) {21};
    \node[above] at (P2r) {21};
    \node[above] at (P3r) {13};
    \node[above] at (P4) {32};

    \node[below] at (Q1) {12};
    \node[below] at (Q2) {31};
    \node[below] at (Q3) {23};
    \node[below] at (Q3r) {23};
    \node[below] at (Q4) {12};

    \node[above=16pt] at ($(P1)!0.5!(P2)$) {$\delta$};
    \draw[->] ($(P1)+(0,0.22)$) -- ($(P2)+(0,0.22)$);

    \node[below=16pt] at ($(Q1)!0.5!(Q3)$) {$\gamma$};
    \draw[<-] ($(Q1)-(0,0.22)$) -- ($(Q3)-(0,0.22)$);

    \node[above=16pt] at ($(P2r)!0.5!(P4)$) {$\gamma$};
    \draw[->] ($(P2r)+(0,0.22)$) -- ($(P4)+(0,0.22)$);

    \node[below=16pt] at ($(Q3r)!0.5!(Q4)$) {$\delta$};
    \draw[<-] ($(Q3r)-(0,0.22)$) -- ($(Q4)-(0,0.22)$);
  \end{tikzpicture}

  \caption{A fundamental polygon equivalent to that shown in Figure \ref{fig:nerve-fundamental-polygon} cut along $[21,23]$ with oppositely orientated cycles.}
  \label{fig:nerve-fundamental-polygon-4}
\end{figure}

\begin{figure}[H]
  \centering
  \begin{tikzpicture}[scale=2.5]
    \coordinate (P1) at (0,1);
    \coordinate (P2) at (0.333,1);
    \coordinate (P3) at (0.667,1);

    \coordinate (Q1) at (0,0);
    \coordinate (Q2) at (0.333,0);
    \coordinate (Q3) at (0.667,0);

    \coordinate (shift) at (0.6,0);

    \coordinate (P2r) at ($(P2)+(shift)$);
    \coordinate (P3r) at ($(P3)+(shift)$);
    \coordinate (P4)  at ($(1,1)+(shift)$);

    \coordinate (Q3r) at ($(Q3)+(shift)$);
    \coordinate (Q4)  at ($(1,0)+(shift)$);

    \draw (Q1) -- (P1) -- (P2) -- (Q3) -- (Q2) -- cycle;

    \draw (Q3r) -- (Q4) -- (P4) -- (P3r) -- (P2r) -- cycle;

    \node[left]  at ($(Q1)!0.5!(P1)$) {$\alpha$};
    \node[right] at ($(Q4)!0.5!(P4)$) {$\alpha$};

    \draw[thick] (P2) -- (Q3);
    \draw[thick] (P2r) -- (Q3r);

    \node[right]  at ($(P2)!0.5!(Q3)$) {$\beta$};
    \node[left] at ($(P2r)!0.5!(Q3r)$) {$\beta$};



    \foreach \pt in {P2,Q2}{
        \fill (\pt) circle (0.01);
      }
    \foreach \pt in {P1,P2r,P3r,P4,Q1,Q3,Q3r,Q4}{
        \fill (\pt) circle (0.02);
      }

    \node[above] at (P1) {21};
    \node[above] at (P2) {13};
    \node[above] at (P2r) {13};
    \node[above] at (P3r) {32};
    \node[above] at (P4) {21};

    \node[below] at (Q1) {31};
    \node[below] at (Q2) {23};
    \node[below] at (Q3) {12};
    \node[below] at (Q3r) {12};
    \node[below] at (Q4) {31};

    \node[above=16pt] at ($(P1)!0.5!(P2)$) {$\delta$};
    \draw[->] ($(P1)+(0,0.22)$) -- ($(P2)+(0,0.22)$);

    \node[below=16pt] at ($(Q1)!0.5!(Q3)$) {$\gamma$};
    \draw[<-] ($(Q1)-(0,0.22)$) -- ($(Q3)-(0,0.22)$);

    \node[above=16pt] at ($(P2r)!0.5!(P4)$) {$\gamma$};
    \draw[->] ($(P2r)+(0,0.22)$) -- ($(P4)+(0,0.22)$);

    \node[below=16pt] at ($(Q3r)!0.5!(Q4)$) {$\delta$};
    \draw[<-] ($(Q3r)-(0,0.22)$) -- ($(Q4)-(0,0.22)$);
  \end{tikzpicture}

  \caption{A fundamental polygon equivalent to that shown in Figure \ref{fig:nerve-fundamental-polygon} cut along $[13,12]$ with oppositely orientated cycles.}
  \label{fig:nerve-fundamental-polygon-5}
\end{figure}

This completes our proof of the validity of oppositely orienting the subspaces representing preference cycles before identification.

\subsection[\appendixname~\thesection]{Proof of Theorem~\ref{theorem:realised-contradictory}}\label{appendix:topology}

\renewcommand{\thesection}{\Alph{section}}

In this section, we complete the proof of Theorem~\ref{theorem:realised-contradictory} by explicitly showing that the identification of the oppositely oriented disks $\DTop$ and $\DBot$ in the capped cylinder $\overline{C} \coloneqq \DTop \cup C \cup \DBot$ produce a surface $X$ that is homeomorphic to the real projective plane ($\RP$). We prove this by first proving that $X$ is a closed, non-orientable surface, and then, we use the classification theorem for closed surfaces to conclude that the non-orientable surface $X$ must be homeomorphic to $\RP$.

To prove that $X$ is a surface and non-orientable, we use the well-known approach of constructing a $2$-sheeted covering map $q: \XCover \rightarrow X$ from a connected, orientable surface $\XCover$. To summarise, a $2$-sheeted covering map $q$ is a $2$-to-$1$ covering map, where the pre-image of each point in $X$ has a point from each orientation of $\XCover$. Firstly this implies that $X$ is a surface, since $q$ is a covering map from a surface $\XCover$, so it follows that $X$ is locally homeomorphic to $\mathbb{R}^2$ and hence is a surface. Secondly $X$ is non-orientable, since given a 2-sheeted covering map $q$, it is well-known that $X$ is non-orientable if and only if $\XCover$ is connected.

\begin{Definition}\label{definition:covering-map}
  Let $X$ and $Y$ be topological spaces. A continuous surjective map $p : X \to Y$ is called a \textbf{covering map} if, for every point $y \in Y$, there exists an open neighbourhood $U \subseteq Y$ of $y$ such that $p^{-1}(U)$ is a disjoint union of open sets in $X$, each mapped homeomorphically onto $U$ by $p$.
\end{Definition}

\begin{Definition}\label{definition:2-sheeted-covering-map}
  Let $S$ be a connected surface (2-manifold). A \textbf{2-sheeted covering map} of $S$ is a pair $(\SCover, p)$ where the following is the case:
  \begin{enumerate}
    \item $\SCover$ is an orientable surface.
    \item $p: \SCover\to S$ is a continuous, surjective map.
    \item $p$ is a covering map such that each fiber $p^{-1}(x)$ consists of exactly two points.
    \item The covering map $p$ is \term{regular}, meaning that there exists a map $\tau: \SCover \rightarrow \SCover$ called a \term{deck transformation} with $p \circ \tau = p$ that exchanges points per fiber and is orientation-reversing.
  \end{enumerate}
\end{Definition}

\begin{Theorem}\label{theorem:two-sheeted-non-orientability}
  Let $S$ be a surface with a 2-sheeted covering map $(\SCover, p)$. $S$ is non-orientable if and only if $\SCover$ is connected.
\end{Theorem}

See~\cite{orientations-textbook} [Theorem 15.4] for a proof of this well-known result. {Examples of $2$-sheeted covering maps from connected surfaces to well-known non-orientable surfaces include maps from cylinders to Möbius strips, from tori to Klein bottles, and from spheres to real projective planes.}

We now proceed to formally define the construction of Theorem~\ref{theorem:realised-contradictory}, {and we show that the constructed surface is homeomorphic to the real projective plane ($\RP$).}

\begin{Construction}\label{construction}
  Recall that the construction underlying Theorem~\ref{theorem:realised-contradictory} is performed as follows. We begin with a cylinder $C = S^1\times [0,1]$ and attach disks $\DTop$ and $\DBot$ to its top and bottom boundaries, which results in the \textit{capped cylinder} $\overline{C} = \DTop \cup C \cup \DBot$. Next, we identify $\DTop$ and $\DBot$ in $\overline{C}$ with an orientation-reversing homeomorphism along $\DTop$ and $\DBot$. We denote the resulting surface $X$, formally defined as the image of a continuous map $q$ on $\overline{C}$ satisfying the following:
  \begin{enumerate}
    \item $q$ is the identity map on $\DTop$, and $q$ identifies $\DBot$ and $\DTop$, i.e., $q(\DBot) = q(\DTop) = \DTop$
    \item $q$ identifies the boundaries of the two disks, i.e., $q(\partial \DBot) = q(\partial \DTop) = \partial \DTop$
    \item $q$ is orientation-reversing on $\partial \DBot$.
  \end{enumerate}

  We proceed to define $q$ in two parts: firstly, a mapping on the \textit{uncapped cylinder} $C = S^1 \times [0,1]$ and then a mapping on the caps $\DTop \cup \DBot$. In order for this definition of $q$ to be well-defined, we must additionally show that these two mappings agree on the shared boundary $\partial(\DTop \cup \DBot) = \partial C$.

  To begin, we parametrise $C$ by coordinates $(\theta,t) \in [0, 2\pi] \times [-\tfrac{1}{2},\tfrac{1}{2}]$ in a standard manner. Intuitively, $t$ denotes how \say{far up} the cylinder the coordinate is, and $\theta$ denotes where around the circumference of the cylinder the coordinate is. This allows us to define $q$ on $C$ as follows.
  \begin{equation}\label{equation:decomposition-1}
    q(\theta,t) = (\ -\theta \bmod 2\pi,\ |t|\ ).
  \end{equation}
  This map is clearly continuous and on the boundary $\partial C = \{(\theta,b)\mid \theta\in[0,2\pi],\ b\in\{-\tfrac{1}{2},\tfrac{1}{2}\}\}$; $q$ clearly identifies the circle at $b = -\tfrac{1}{2}$ with the circle at $b = \tfrac{1}{2}$ with an orientation-reversing twist.

  Next, we parametrise the disks $\DTop \cup \DBot$ by coordinates $(\varphi,r,b) \in [0, 2\pi] \times [0,1] \times \{\text{bot},\text{top}\}$. Intuitively, $r$ denotes where the coordinate lies between the origin and boundary of the disk, $\varphi$ denotes where around the circumference of the disk the coordinate radiates to, and $b$ denotes whether the coordinate is in $\DBot$ or $\DTop$.

  \noindent This allows us to define $q$ on $\DTop \cup \DBot$ as follows.
  \begin{equation}\label{equation:decomposition-2}
    q(\varphi,r,b) = \begin{cases*}
      (\varphi,r,\text{top})             & \text{when $b = \text{top}$}             \\
      (-\varphi \bmod 2\pi,r,\text{top}) & \text{otherwise, i.e., $b = \text{bot}$}
    \end{cases*}
  \end{equation}
  This map is likewise continuous and identifies the disks in an orientation-reversing manner. Finally, we verify that Equations (\ref{equation:decomposition-1}) and (\ref{equation:decomposition-2}) agree on their boundaries post-identification; this follows from the equivalence of $(\varphi, 1, \text{top}) \sim (\theta, \tfrac{1}{2})$ when $\varphi = \theta$ in accordance with $\partial(\DTop \cup \DBot) = \partial C$.
\end{Construction}

\pagebreak

\begin{Theorem}\label{theorem:main-construction}
  Construction~\ref{construction} defines a surface $X$ that is homeomorphic to $\RP$.
\end{Theorem}
\begin{proof}
  We first prove that $X$ is non-orientable by defining a 2-sheeted covering map $q: \XCover \rightarrow X$ from a connected, orientable surface $\XCover$. Indeed, we take $\XCover = \overline{C}$, and we proceed to show that the function $q: \XCover \rightarrow X$ as defined in Construction~\ref{construction} is a 2-sheeted covering map. The conditions of Definition~\ref{definition:2-sheeted-covering-map} are verified as follows.

  Firstly, the capped cylinder $\XCover$ is clearly orientable, as it is homeomorphic to the sphere $S^2$.

  Secondly, $q$ is continuous, as it is defined by the continuous mappings of Equations (\ref{equation:decomposition-1}) and (\ref{equation:decomposition-2}), and it is surjective by definition of $X$ being defined as the image of those mappings.

  Finally, it remains to show that $q$ is a $2$-to-$1$ covering map such that the map defined by interchanging the two points of every fiber comprises an orientation-reversing deck transformation. Indeed, recall that $q$ is defined by $2$-to-$1$ identifications of the form $(\theta, b) \leftrightarrow (-\theta \bmod 2\pi, -b)$ or $(\varphi, r, \text{bot}) \leftrightarrow (-\varphi \bmod 2\pi, r, \text{top})$. One can construct a ($1$-to-$1$) homeomorphism $\XCover \cong S^2$ such that the fibres of $q$ correspond to antipodal identifications $(x,y,z) \leftrightarrow (-x, -y, -z)$ on $S^2$. Visually, this corresponds to splitting a sphere into three regions: an equatorial band corresponding to the cylinder $C \subset \overline{C}$ and \textit{northern} and \textit{southern} hemispheres corresponding to $\DTop$ and $\DBot$. Thus, interchanging the two points of every fiber is orientation-reversing and hence a deck transformation. Furthermore, $q$ is a covering map, as for every point $x \in X$, one can construct an open neighbourhood $U$ for which its preimage under $q$ consists of two \textit{equally sized}, disjoint, open neighbourhoods in $\overline{C}$ that surround the two \say{antipodes} comprising $q^{-1}(x)$ (see Figure~\ref{fig:appendix}).

  \begin{figure}[H]
    \centering
    \includegraphics[width=0.45\textwidth]{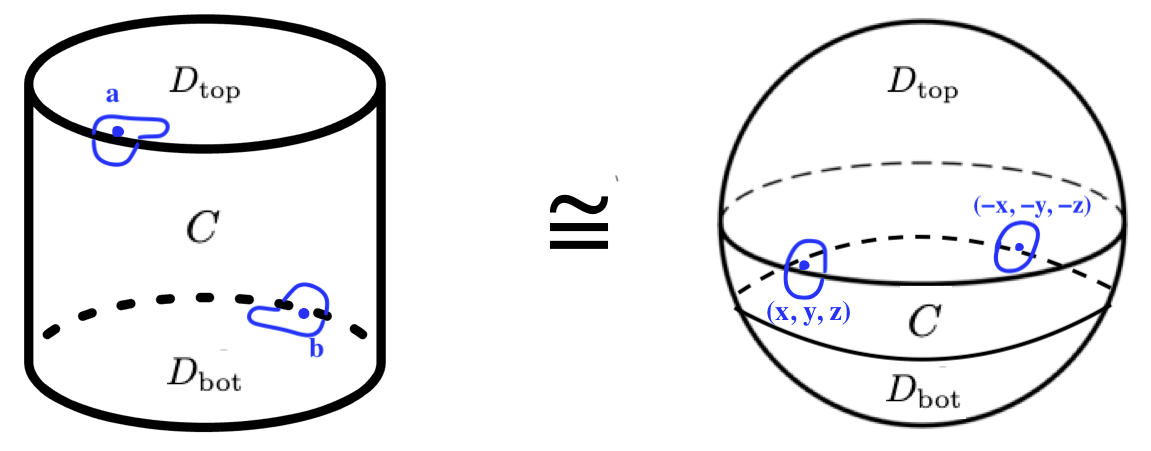}
    \caption{Correspondence of a fiber $\{a,b\}$ identified by $q$ and antipodal points $(x,y,z)$ and $(-x,-y,-z)$ on a sphere (irrespective of whether $a$ and $b$ are both in $\DTop \cup \DBot$ or both in $C$).}
    \label{fig:appendix}
  \end{figure}

  Thus, $q$ is a $2$-sheeted covering map. We can now conclude that $X$ is a surface and non-orientable. Firstly, $X$ is a surface, since $q$ is a covering map from a surface $\XCover$, so it follows that $X$ is locally homeomorphic to $\mathbb{R}^2$ and hence is a surface; $X$ is also closed because $\XCover$ is. Secondly, $X$ is non-orientable since $\XCover = \overline{C} \cong S^2$ is connected and orientable, so by Theorem~\ref{theorem:two-sheeted-non-orientability}, it follows that $X$ is connected and non-orientable. Finally, since $q$ is $2$-to-$1$, the Euler characteristic $\chi(X)$ of $X$ must then satisfy the following relation.
  \begin{equation*}
    \chi(X) = \tfrac{1}{2}\chi(\XCover) = \tfrac{1}{2}\chi(S^2) = 1.
  \end{equation*}
  By the classification theorem for closed surfaces, the unique connected, non-orientable surface with Euler characteristic $1$ is $\RP$. Hence, $X$ is homeomorphic to $\RP$.\\
\end{proof}

\subsection[\appendixname~\thesection]{\text{Arrow's Impossibility Theorem as Non-Orientability} \linebreak ($>$3 Alternative Case)\label{appendix:more-alternatives}}

In this section, we prove the general cases of Theorems~\ref{theorem:arrow-reduction-1} and~\ref{theorem:stronger-equivalence-1}, i.e., for 3 or more alternatives and allowing indifference between alternatives too. We begin by defining SWFCs for 3 or more alternatives with indifference, as well as the corresponding Unrestricted Domain, Unanimity, IIA and Non-Dictatorship conditions. These definitions are straightforward reformulations of those found in~\cite[Section 3 and Appendix A]{paper0-arxiv};
omitting the indifference symbol $\indiff$, one recovers D'Antoni's definitions from~\cite[Section 3]{dantoni}.

To begin, fix a set $\A = \{1,2,3,\dots,a\}$ of 3 or more alternatives;  let $\P{\A}$ be the set of weak orders on $\A$, and let $\PN{\A}$ be the set of $N$-individual profiles on $\A$ for $N \geq 2$ individuals. Then, let $\WeakPrefs{\A}$ be the set of $\binom{|\A|}{2}$-length tuples with values in $\{0,\indiff,1\}$, mirroring Section~\ref{subsection:generalised-social-welfare-functions} with the addition of $\indiff$ to denote indifference (see~\cite[Definition A.1]{paper0-arxiv}.)

For $A = \binom{|\A|}{2}$, an $N$ individual profile on $\A$ can then be represented by a ternary-valued $A \times N$ matrix, where the rows of these matrices record each individual’s preferences on a single pair of alternatives, and each column records an individual’s preferences~\cite[p.4]{paper0-arxiv}.

Recall that, in Section~\ref{subsection:reference-orientation}, we introduced a quotienting (i.e., identification) process for the contradictory preference cycles of $\WeakPrefs{\3}$ to form $\PCont{\3} \coloneqq \P{\3} \sqcup\{{\textbf{c}}\}$. Moreover, recall that a \term{Social Welfare Function (SWF)} in the sense of~\cite{dantoni,paper0-arxiv} is any function of the form $\PPN{\A} \rightarrow \Prefs{\A}$ in the strict case and of the form $\PN{\A} \rightarrow \WeakPrefs{\A}$ in the general case. Furthermore, for 3-alternative sets $\3$, we further define a \term{Social Welfare Function with Contradictory Cycles (SWFC)} as any function of the form $\PPN{\3} \rightarrow \PPCont{\3}$ in the strict case and $\PN{\3} \rightarrow \PCont{\3}$ in the general case.

Because we will solely work with SWFs/SWFCs that satisfy IIA, we can simplify several definitions by first defining IIA, and then, we define every other property under the assumption of IIA thereafter.

\parsplit

\begin{Definition}[IIA]\label{definition:iia}
  Let $w: \PN{\A} \rightarrow \WeakPrefs{\A}$ be an SWF: $w$ satisfies \term{IIA} if $w$ is the product of some series of functions $s_1, s_2, \dots s_A: \{0,\indiff,1\}^N \rightarrow \{0,\indiff,1\}$ for $A = \binom{|\A|}{2}$. That is, for any $p \in \PN{\A}$ for which its matrix representation has rows $u_1,u_2,\dots,u_A$, we have that $w(p) = (s_1(u_1),s_2(u_2),\dots,s_A(u_A))$.
\end{Definition}

\parsplit

\begin{Definition}[Unanimity and Non-Dictatorship]\label{definition:unanimity-nd}
  Let $w: \PN{\A} \rightarrow \WeakPrefs{\A}$ be an SWF satisfying IIA with respect to the series of functions $s_1, s_2, \dots s_A$. The SWF $w$ satisfies the following:
  \begin{itemize}
    \item\textbf{Unanimity:} If $\forall\ j \in \{1,2,\dots,A\}$ and $\forall\ x \in \{0,1\}$, when denoting $\Delta x = \underbrace{(x,x,\dots, x)}_{N \text{ times}}$, we have that $s_j(\Delta x) = x$.
    \item\textbf{Dictatorship at Individual i:} If $\forall j \in \{1,2,\dots,A\}$ and $\forall\ (u_1,\dots,u_i,\dots,u_N) \in \{0,\indiff,1\}^N$ with $u_i \in \{0,1\}$, we have that $s_j(u_1,\dots,u_i,\dots,u_N) = u_i$. If no such individual $i$ exists, we say that $w$ satisfies \term{Non-Dictatorship}.
  \end{itemize}
\end{Definition}

\noindent For any subset of 3 alternatives, $\3 \subseteq \A$, there are \textit{restrictions} $r_{\3}: \PN{\A} \rightarrow \PN{\3}$ and, likewise, $s_{\3}: \WeakPrefs{\A} \rightarrow \WeakPrefs{\3}$ that \textit{drop} all pairwise comparisons in $\A$ outside $\3$.

\parsplit

\begin{Example}\label{example:restriction}
  Let $\A =\{a,b,c,d\}$ be a set of 4 alternatives, and we consider profiles of 2 individuals. If $p$ is the profile $(a \prec b \prec c \prec d,\ d \prec b \prec a \prec c)$, then the restriction of $r_{\3}(p)$ for $\3 = \{a,b,d\}$ would be $(a \prec b \prec d,\ d \prec b \prec a)$. Likewise, a preference cycle $r$ on $\A$ given by $a \prec b \prec c \prec d \prec a$ has restriction $s_{\3}(r)$ given by $a \prec b \prec d \prec a$.
\end{Example}

\parsplit

Hence, IIA implies that there exists a Social Welfare Function $w_{\3}: \PN{\3} \rightarrow \WeakPrefs{\3}$ such that the following diagram commutes.

\begin{equation}
  \begin{tikzcd}[row sep=3em, column sep=4em]
    \PN{\A} \arrow[d,"w"'] \arrow[r,"r_{\3}"] & \PN{\3} \arrow[d,"w_{\3}"] \\
    \WeakPrefs{\A} \arrow[r,"s_{\3}"'] & \WeakPrefs{\3}
  \end{tikzcd}\label{equation:pre-commutative}
\end{equation}

\parsplit

\begin{Definition}[Unrestricted Domain]\label{definition:ud}
  Let $w: \PN{\A} \rightarrow \WeakPrefs{\A}$ be an SWF: $w$ satisfies the \term{Unrestricted Domain} if $\forall p \in \PN{\A}$ we have that, for every 3-alternative subset, $\3 \subseteq \A$: $s_{\3}(w(p))$ is not a preference cycle.
\end{Definition}

\parsplit

Now, let $w: \PN{\A} \rightarrow \WeakPrefs{\A}$ be an SWF that only satisfies IIA and Unanimity and thus potentially aggregates some profiles to relations that contain preference cycles. The quotienting of $\Prefs{\3}$ to $\PPCont{\3}$ defines a surjection $\sigma: \WeakPrefs{\3} \rightarrow \PCont{\3}$, which is the identity on weak orders, and it maps the preference cycles to \textbf{c}. Hence, we redefine $w_{\3}$ and $s_{\3}$ in Equation (\ref{equation:pre-commutative}) by post-composing with $\sigma$ to produce the following.

\begin{equation}
  \begin{tikzcd}[row sep=3em, column sep=4em]
    \PN{\A} \arrow[d,"w"'] \arrow[r,"r_{\3}"] & \PN{\3} \arrow[d,"w_{\3}"] \\
    \WeakPrefs{\A} \arrow[r,"s_{\3}"'] & \PCont{\3}
  \end{tikzcd}\label{equation:commutative}
\end{equation}

Next, we restrict $w$ to profiles $\D \subseteq \PN{\A}$, for which their valid aggregations are strict on $\3$. Then, we denote $\psi: \D \rightarrow \WeakPrefs{\A}$ as the restriction of $w$ to $\D$ and $\psi_{\3}$ as the restriction of $w_{\3}$ to $D(\3) \coloneqq \{r_{\3}(p) \mid p \in \D\}$ so that the following diagram commutes.
\begin{equation}
  \begin{tikzcd}[row sep=3em, column sep=4em]
    \D \arrow[d,"\psi"'] \arrow[r,"r_{\3}"] & \D(\3) \arrow[d,"\psi_{\3}"] \\
    im{(\psi)} \subseteq \WeakPrefs{\A} \arrow[r,"s_{\3}"'] & \PPCont{\3}
  \end{tikzcd}\label{equation:commutative-2}
\end{equation}
In other words, $\psi$ is the restriction of $w$ to all profiles whose aggregates are strict orders or relations containing a preference cycle, and hence $s_{\3}$ restricted to $im(\psi)$ has codomain $\PPCont{\3}$ rather than $\PCont{\3}$.

Similarly to Lemma~\ref{lemma:two-options}, we have that $im(\psi_{\3}) = \PP{\3}$ or $im(\psi_{\3}) = \PPCont{\3}$ because, by the Unanimity property of $w$, we have that $x \in \PP{\3}$: $w(\Delta x) = x$, which implies that $\Delta x \in \D(\3)$, further implying that $x \in im(\psi_{\3})$. This allows us to define an Arrovian Topological Model of $\psi_{\3}$ and $\psi_{\3}$ with Orientable Social Choices, just as we did in Definitions~\ref{definition:arrovian-topological-model} and \ref{definition:orientable-social-choices} as follows.

\parsplit

\begin{Definition}\label{definition:psi-versions}
  We define an \term{Arrovian Topological Model} of $\psi_{\3}$ as any surface that is a topological model of $\PPCont{\3}$ with realised preference cycles when $\psi_{\3} = \PPCont{\3}$ and a topological model of $\PP{\3}$ otherwise, i.e., when $im(\psi_{\3}) = \PP{\3}$. Furthermore, if all Arrovian Topological Models of $\psi_{\3}$ are orientable, we say that $\psi_{\3}$ has \term{Orientable Social Choices} and \term{Non-Orientable Social Choices} otherwise.
\end{Definition}

\parsplit

\begin{Proposition}\label{proposition:equivalent-conditions-2}
  Let $w: \PN{\A} \rightarrow \WeakPrefs{\A}$ be an SWF on a set $\A$ with 3 or more alternatives and $N \geq 2$ individuals. The SWF $w$ satisfies the Unrestricted Domain if and only if for every subset of 3 alternatives, $\3 \subseteq \A$: $\psi_{\3}$ has Orientable Social Choices.
\end{Proposition}
\begin{proof}
  ($\implies$) If $w$ satisfies the Unrestricted Domain, then, for every $\3$, $im(\psi_{\3}) = \PP{\3}$. So, arguing as in Theorem~\ref{proposition:equivalent-conditions}, we have that $\psi_{\3}$ has Orientable Social Choices because all its Arrovian Topological Models are homeomorphic to the annulus.

  ($\impliedby$) Assume by way of contradiction that, for every $\3 \subseteq \A$: $\psi_{\3}$, has Orientable Social Choices but $w$ does not satisfy the Unrestricted Domain. This implies that there exists a 3-alternative subset $\threeast \subseteq \A$ such that $im(\psi_{\threeast}) = \PPCont{\threeast}$. Hence, arguing as in Theorem~\ref{proposition:equivalent-conditions}, we have that $\psi_{\threeast}$ does not have Orientable Social Choices because its Arrovian Topological Models must be homeomorphic to the real projective plane, which is a contradiction.\\
\end{proof}

\parsplit

\begin{Theorem}[Reformulation of Arrow's Impossibility Theorem as Non-Orientability]\label{theorem:arrow-reduction-2}~\\
  The following two statements are equivalent:
  \begin{enumerate}
    \item No SWF $w: \PN{\A} \rightarrow \WeakPrefs{\A}$ satisfies all of Unanimity, IIA, Non-Dictatorship, and the Unrestricted Domain.
    \item No SWF $w: \PN{\A} \rightarrow \WeakPrefs{\A}$ satisfies all of Unanimity, IIA, and Non-Dictatorship, and for every subset of 3 alternatives $\3 \subseteq \A$, $\psi_{\3}$ has Orientable Social Choices.
  \end{enumerate}
\end{Theorem}
\begin{proof}
  We simply apply Proposition~\ref{proposition:equivalent-conditions-2} to interchange the conditions of $w$ satisfying the Unrestricted Domain and every $\psi_{\3}$ having Orientable Social Choices.\\
\end{proof}

The following results establish an equivalence between the conditions of Non-Dictatorship and Non-Orientable Social Choices, but only for social welfare functions without indifference. This limitation reflects that a dictatorial social welfare function may still produce preference cycles for profiles where the dictator has no strict preferences.

\parsplit

\begin{Lemma}\label{lemma:dictator-to-ud}
  Let $w: \PPN{\A} \rightarrow \Prefs{\A}$ be an SWF satisfying IIA. If $w$ has a dictator at $i$ for some $i \in \{1,2,\dots, N\}$ then $\forall p \in \PPN{\A}$ we have that $w(p)$ cannot have a preference cycle.
\end{Lemma}
\begin{proof}
  This simply follows from the fact that if $i$'s preferences are given by a strict order $x \in \PP{\A}$, we have that $w$ having a dictator at $i$ implies that $w(p) = x$, and strict orders of course have no preference cycles.\\
\end{proof}

\begin{Theorem}[Non-Dictatorship as Non-Orientability]\label{theorem:stronger-equivalence-2}
  Let $w: \PPN{\A} \rightarrow \Prefs{\A}$ be an SWF that satisfies IIA and Unanimity. The SWF $w$ additionally satisfies Non-Dictatorship if and only if there exists a subset of 3 alternatives $\threeast \subseteq \A$ such that $\psi_{\threeast}$ has Non-Orientable Social Choices.
\end{Theorem}
\begin{proof}
  ($\implies$) The SWF $w$ satisfies Non-Dictatorship implies there exists a subset of 3 alternatives $\threeast \subseteq \A$ such that $im(\psi_{\threeast}) = \PPCont{\3}$, which implies that it has Non-Orientable Social Choices.

  ($\impliedby$) If there is some $\threeast \subseteq \A$ such that $im(\psi_{\threeast}) = \PPCont{\3}$ has Non-Orientable Social Choices, then by Theorem~\ref{theorem:arrow-reduction-2}, $w$ does not satisfy the Unrestricted Domain, which by the contrapositive to Lemma~\ref{lemma:dictator-to-ud} implies that $w$ must satisfy Non-Dictatorship.\\
\end{proof}

\end{document}